\documentclass[11pt]{article}

\usepackage{latexsym,mathrsfs}
\usepackage{amsmath,amssymb}
\usepackage{amsthm,enumerate,verbatim}
\usepackage{amsfonts}
\usepackage{graphicx}
\usepackage{algorithm}
\usepackage{algorithmic}
\usepackage{url}

\newtheorem*{definition*}{Definition}

\newtheorem{remark}{Remark}

\DeclareMathOperator{\argmin}{argmin}

\title{ Accelerating Nonnegative Matrix Factorization Algorithms \\ using Extrapolation\thanks{The authors acknowledge the support of the European Research Council (ERC starting grant n$^\text{o}$ 679515).} }
\date{}

\author{Andersen Man Shun Ang$^{\dagger}$ \qquad Nicolas Gillis\thanks{E-mails:  \{manshun.ang, nicolas.gillis\}@umons.ac.be.} \\   
Department of Mathematics and Operational Research \\ 
Facult\'e Polytechnique, Universit\'e de Mons \\ 
Rue de Houdain 9, 7000 Mons, Belgium
}

\begin{document}

\maketitle

\begin{abstract} 
In this paper, we propose a general framework to accelerate significantly the algorithms for nonnegative matrix factorization (NMF). This framework is inspired from the extrapolation scheme used to accelerate gradient methods in convex optimization and from the method of parallel tangents. However, the use of extrapolation in the context of the exact two-block coordinate descent algorithms tackling the non-convex NMF problems is novel. We illustrate the performance of this approach on two state-of-the-art NMF algorithms, namely, accelerated hierarchical alternating least squares (A-HALS) and alternating nonnegative least squares (ANLS), using synthetic, image and document data sets.  
\end{abstract} 


\section{Introduction} 
Given an input data matrix $X \in \mathbb{R}^{m \times n}$ and a factorization rank $r$, we consider in this paper the following optimization problem 
\begin{equation} \label{nmf}
\min_{W \in \mathbb{R}^{m \times r}, H \in \mathbb{R}^{r \times n}} \; ||X - WH||_F^2 
\quad \text{ such that } \quad 
W \geq 0 \text{ and } H \geq 0. 
\end{equation}
This problem is referred to as nonnegative matrix factorization (NMF) and has been shown to be useful in many applications such as image analysis and documents classification~\cite{lee1999learning}. 
Note that there exists many variants of~\eqref{nmf} using other objective functions and additional contraints and/or penalty term on $W$ and $H$; see~\cite{cichocki2009nonnegative, gillis2014, gillis2017, fu2018nonnegative} and the references therein for more details about NMF models and their applications.

\paragraph{Algorithms for NMF} The focus of this paper is algorithm design for~\eqref{nmf}. Almost all algorithms for NMF use a two-block coordinate descent scheme by optimizing alternatively over $W$ for $H$ fixed and vice versa; see Algorithm~\ref{nmfalgo}. 
\algsetup{indent=2em}
\begin{algorithm}[ht!]
\caption{Framework for most NMF algorithms} \label{nmfalgo}
\begin{algorithmic}[1]
\REQUIRE
%The (non-convex) function $f(x)$,
%the feasible set $\mathcal{X}$,
An input matrix $X \in \mathbb{R}^{m \times n}$, 
an initialization 
$W \in \mathbb{R}^{m \times r}_+$, 
$H \in \mathbb{R}^{m \times r}_+$. 

\ENSURE An approximate solution $(W,H)$ to NMF.  \medskip

\FOR{$k = 1, 2, \dots$}
\STATE Update $W$ using a NNLS algorithm to minimize 
%\[
$||X - W H||_F^2$ with $W \geq 0$. 
%\]
\STATE Update $H$ using an NNLS algorithm to minimize 
%\[
$||X - W H||_F^2$  with $H \geq 0$. 
%\] 
\ENDFOR
\end{algorithmic}
\end{algorithm} 
By symmetry, since $||X-WH||_F = ||X^T-H^TW^T||_F$, the updates of $W$ and $H$ are usually based on the same strategy. Looking at the subproblem for $H$, the following nonnegative least squares (NNLS) problem 
\begin{equation} \label{nnls}
\min_{H \geq 0} \; ||X - WH||_F^2 
\end{equation}
needs to be solved exactly or approximately. 
The most popular approaches in the NMF community to solve~\eqref{nnls} are 
multiplicative updates~\cite{lee1999learning}, 
active-set methods that solve \eqref{nnls} exactly~\cite{kim2008nonnegative, kim2011fast}, 
projected gradient methods~\cite{lin2007projected, nenmf}, 
and exact block coordinate descent (BCD) methods~\cite{Cic, Cic4, hsieh2011fast, gillis2012accelerated, chow2017cyclic}. 
Among these approaches, exact BCD schemes have been shown to be the most effective in most situations~\cite{kim2014algorithms}. The reason is that the optimal update of a single row of $H$, the others being fixed, admits a simple closed-form solution: we have for all $k$ that 
\begin{align*}
\argmin_{H(k,:) \geq 0} ||X - WH||_F^2 
&  = 
\max \left( 0 , \frac{ W(:,k)^T ( X - \sum_{j \neq k} W(:,j) H(j,:) ) }{||W(:,k)||_2^2} \right) \\ 
&  = 
\max \left( 0 , \frac{  W(:,k)^T X - \sum_{j \neq k} \left( W(:,k)^T W(:,j) \right) H(j,:) }{||W(:,k)||_2^2} \right).  
\end{align*} 
The algorithm using these updates is referred to as hierarchical alternating least squares (HALS) and  updates the rows of $H$ and the columns of $W$ 
in a sequential way~\cite{Cic, Cic4}.  
HALS has been improved in several ways: 
\begin{itemize} 

\item Selecting the variable to be updated in order to reduce the objective function the most (Gauss-Seidel coordinate descent)~\cite{hsieh2011fast}. 

\item Updating the rows of $H$ several times before updating $W$ (and similarly for the columns of $W$) as the computation of $W^TW$ and $W^TX$ can be reused which allows a significant acceleration of HALS~\cite{gillis2012accelerated}. This variant is referred to as accelerated HALS (A-HALS). 

\item  Using random shuffling instead of the cyclic updates of the rows of $H$ which leads in general to better performances~\cite{chow2017cyclic}. However, when combined with the above strategies to accelerate HALS, we have observed that the improvement is negligible.

\end{itemize} 
More recently, HALS was also accelerated using randomized sampling techniques~\cite{erichson2018randomized}. 

Although the acceleration scheme proposed in this paper can potentially be applied to any NMF algorithm, we will focus for simplicity on two algorithms: 
\begin{enumerate}
\item A-HALS which is, as explained above, arguably one of the most efficient NMF algorithm, and 

\item Alternating nonnegative least squares (ANLS) which is Algorithm~\ref{nmfalgo} where the NNLS subproblems~\eqref{nnls} are solved exactly. To solve the NNLS subproblems, 
we use the active-set method from~\cite{kim2011fast} which is one of the most efficient strategy for NNLS~\cite{kim2014algorithms}. 

\end{enumerate}

\paragraph{Outline of the paper} 
In this paper, we introduce a general framework to accelerate NMF algorithms. 
This framework, described in Sections~\ref{accel} and~\ref{extranmf}, is closely related to the extrapolation scheme usually used in the context of gradient descent methods. 
We use it here in the context of exact BCD methods applied to NMF.  
The difficulty in using this scheme is in choosing the tuning parameters in the extrapolation, for which we propose a simple strategy in Section~\ref{sec:beta}. We illustrate the effectiveness of this approach on synthetic, image and document data sets in Section~\ref{numexp}.

\section{Acceleration through extrapolation} \label{accel}

\begin{figure}[ht!]
\begin{center}
\includegraphics[width=0.7\textwidth]{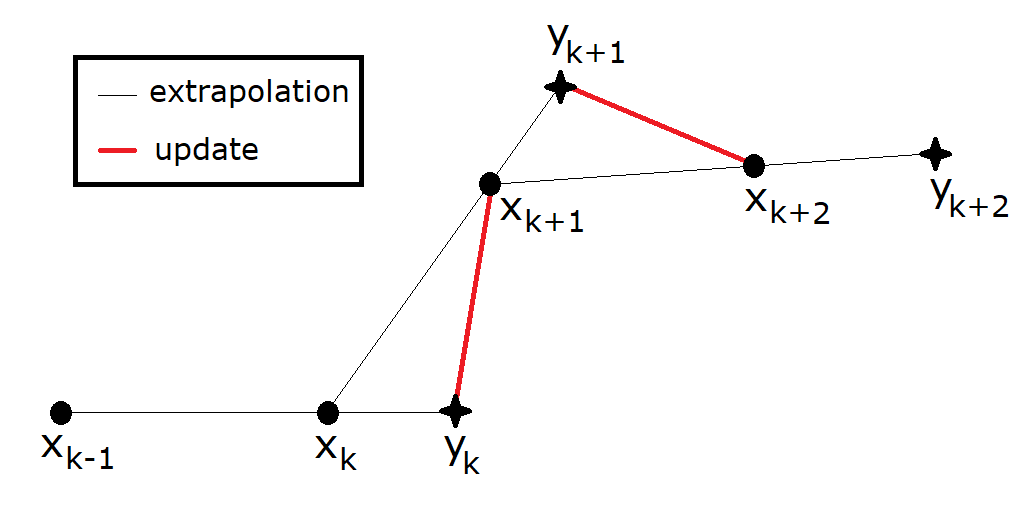} 
\caption{
Illustration of the idea of extrapolation to accelerate optimization schemes.  \label{fig:extrapol}}
\end{center}
\end{figure}

Let us describe the simple extrapolation scheme that we will use to accelerate NMF algorithms. This scheme takes its roots in the so-called method of parallel tangents which is closely related to the conjugate gradient method~\cite[p.~293]{luenberger1984linear}, and the accelerated gradient schemes by Nesterov~\cite{nesterov2013introductory}. 
The idea is the following. 
Let us consider an optimization scheme that computes the next iterate only based on the previous iterate\footnote{Although this assumption is not strictly necessary, it makes more sense otherwise there might be 
a counter effect if the update already takes into account the previous iterates.} 
(e.g., gradient descent or a coordinate descent), 
that is, it updates the $k$th iterate $x_k$ as follows 
\[
x_{k+1} = \text{update}(x_k), 
\]
for some function update(.)  that depend on the objective function and the feasible set.  
For most first-order methods,  
these updates will have a zig-zagging behavior. In particular, gradient descent with exact line search leads to orthogonal search directions~\cite{luenberger1984linear} while search direction of (block) coordinate descent methods are orthogonal by construction. The idea of extrapolation is to define a second sequence of iterates, namely, $y_k$ with $y_0 = x_0$, 
and modify the above scheme as follows 
\[
x_{k+1} = \text{update}(y_k), 
\quad y_{k+1} = x_{k+1}  + \beta_k (x_{k+1} - x_k), 
\]
for some $\beta_k \geq 0$. 
 Note that there are other possibilities for choosing $y_{k+1}$ based on a linear combinations of previous iterates. Figure~\ref{fig:extrapol} illustrates the extrapolation scheme, and allows us to get some intuition: 
 the direction $(x_{k+1} - x_k)$ will be inbetween zigzagging directions obtained with the original update applied to $y_k$'s and will allow to accelerate convergence. For example, we observe on Figure~\ref{fig:extrapol} that the direction $x_{k+2}-x_{k+1}$ is between the directions $x_{k+1}-y_k$ and $x_{k+2}-y_{k+1}$. 
  
 In the case of gradient descent and smooth convex optimization, the above  scheme allows to accelerate convergence of the function values from  $O(1/k)$ to $O(1/k^2)$, and from linear convergence with rate $(1-{\mu/L})$ to rate $(1-\sqrt{\mu/L})$ for strongly convex function with parameter $\mu$ and whose gradient has Lipschitz constant $L$~\cite{nesterov2013introductory}. This scheme has also been used for BCD, and most works focus on the case where the blocks of variables are updated using a gradient or proximal step; see, e.g., \cite{beck2013convergence, xu2013block,  fercoq2015accelerated, chambolle2015remark} and the references therein.  
In the convex case, the $\beta_k$'s can be chosen a priori in order to obtain the theoretical acceleration. However, from a practical point of view, the acceleration will depend on the choice of the $\beta_k$'s which is non-trivial; 
 see, e.g., \cite{o2015adaptive} for a discussion about this issue. 
 %It has to be noted that the iterates $y_k$'s do not necessarily belong to the feasible set. 
Extrapolation has been used more recently in non-convex settings~\cite{xu2013block, o2017behavior, paquette18a}, 
but, as far as we know, not in combination with exact BCD methods. 
Xu and Yin~\cite{xu2013block} used extrapolation in the context of an inexact BCD method where the blocks of variables are updated using a projected gradient method. 
Their approach is different from our as we will use exact BCD.  
Note that Xu and Yin~\cite{xu2013block} applied their technique to NMF which we will compare to ours in Section~\ref{numexp}. 

In the method of parallel tangents, the steps $\beta_k$ are computed using line-search~\cite[p.~293]{luenberger1984linear}. This allows the acceleration scheme to be at least as good as the initial scheme. 
However, as we will see, this is not a good strategy in our case because the optimal $\beta_k$'s will be close to zero (because we use coordinate descent). %We will have to take some `risk' going beyond the optimal $\beta_k$'s to allow acceleration at the next steps. 
 In any case, the choice of the $\beta_k$'s is non-trivial and, as we will see, the acceleration depends on the choice of these parameters. 
 %In fact, even in the convex case, the best choice for the $\beta_k$'s is difficult in practice; see, e.g., \cite{o2015adaptive} for a discussion about this issue. 
 Note that choosing $\beta_k=0$ for all $k$ gives back the original algorithm (no extrapolation), and $\beta_k$ close to one is a very aggressive strategy. \\

 The remainder of this paper is organized as follows: 
 \begin{itemize} 
 
 \item In Section~\ref{extranmf}, we adapt the above extrapolation technique in the context of two-block coordinate NMF algorithms (Algorithm~\ref{nmfalgo}). 
 
 \item In Section~\ref{sec:beta}, we propose a simple strategy for the choice of the parameters $\beta_k$'s. 
  
  \item In Section~\ref{numexp}, we illustrate the acceleration of NMF algorithms on synthetic, image and document data sets. 
 
 \end{itemize}

\section{Extrapolation for NMF algorithms} \label{extranmf}

In this paper, we adapt extrapolation to the two-block coordinate descent strategies of NMF algorithms described in Algorithm~\ref{nmfalgo}. 
Algorithm~\ref{accextra} describes the proposed extrapolation scheme applied to NMF~\eqref{nmf}.  Depending on the choice of the parameter $hp \in \{1,2,3\}$,  %($h$ stands for $H$, $p$ for projection), 
Algorithm~\ref{accextra} corresponds to three different variants of the proposed extrapolation. This is described below through two important questions.

\algsetup{indent=2em}
\begin{algorithm}[ht!]
\caption{Acceleration of Algorithm~\ref{nmfalgo} using extrapolation} \label{accextra}
\begin{algorithmic}[1]
\REQUIRE An input matrix $X \in \mathbb{R}^{m \times n}$, 
an initialization 
$W \in \mathbb{R}^{m \times r}_+$ 
and 
$H \in \mathbb{R}^{m \times r}_+$, 
parameters $hp \in \{1,2,3\}$ (extrapolation/projection of $H$). 
\ENSURE An approximate solution $(W,H)$ to NMF.  \medskip
\STATE $W_y = W$; $H_y = H$; $e(0) = ||X - WH||_F$.
\FOR{$k = 1, 2, \dots$}
\STATE Compute $H_n$ using a NNLS algorithm to minimize $||X - W_y H_n||_F^2$ with $H_n \geq 0$ using $H_y$ as the initial iterate. 
\IF{ $hp \geq 2$ } 
\STATE Extrapolate: $H_y = H_n + \beta_k (H_n - H)$. 
\ENDIF 
\IF{ $hp = 3$ } 
\STATE Project: $H_y = \max\left( 0 , H_y \right)$. 
\ENDIF  
\STATE Compute $W_n$ using a NNLS algorithm to minimize $||X - W_n H_y||_F^2$ with $W_n \geq 0$ using $W_y$ as the initial iterate. 
\STATE Extrapolate: $W_y = W_n + \beta_k  (W_n - W)$. 
\IF{ $hp = 1$ } 
\STATE Extrapolate: $H_y = H_n + \beta_k (H_n - H)$.
\ENDIF 
\STATE Compute the error: e($k$) $= ||X - W_n H_y||_F$. \quad \emph{\% See Remark~\ref{comperr}.}
	\IF{e($k$) $>$ e($k-1$)} 
  		\STATE Restart: $H_y = H$; $W_y = W$. 
  		%\STATE Restart: $\beta = \beta_0$.
  \ELSE 
  		\STATE $H = H_n$; $W = W_n$.
  		%\STATE Update $\beta$; see Section~\ref{sec:beta}. 
	\ENDIF 
\ENDFOR 
\end{algorithmic}
\end{algorithm}

\paragraph{When should we perform the extrapolation?} In case of NMF (and in general for block coordinate descent methods), 
it makes sense to perform the extrapolation scheme after the update of each block of variables so that when we update then next block of variables, the algorithm takes into account the already extrapolated variables; see, e.g.,~\cite{fercoq2015accelerated}. 
However, as we will see in the numerical experiments, this does not necessarily performs best in all cases. This is the first reason why we have added a parameter $hp \in \{1,2,3\}$: 
For $hp = 1$, 
$H$ is extrapolated after the update of $W$, 
otherwise it is extrapolated directly after it has been updated. 
Note that in the former case, the extrapolated matrix $H_y$ is only used as a warm start for the next NNLS update of $H$. For ANLS, it will therefore not play a crucial role since ANLS solves the NNLS subproblem exactly. 
%Note also that, as opposed to Algorithm~\ref{nmfalgo}, choosing $h=0$ makes Algorithm~\ref{accextra} breaks the symmetry in $W$ and $H$. 

\paragraph{Can we guarantee convergence?}  
Under some mild assumptions and/or slight modifications of the algorithm, 
%(or slight modification of the algorithm, using a regulrization ||W-W||_F^2), 
block coordinate descent schemes are guaranteed to converge to stationary points~\cite{hong2017iteration}. Since Algorithm~\ref{accextra} uses extrapolation, we cannot use these results directly. Similarly, we cannot use the result of Xu and Yin~\cite{xu2013block} that use projected gradient steps to update $W$ and $H$. 
 %Algorithm~\ref{accextra} uses extrapolation and does not fall within this class of methods. 
In Algorithm~\ref{accextra}, because $W_y$ and $H_y$ are not necessarily nonnegative, the objective function is not  guaranteed to decrease at each step. In fact, step~16 of  Algorithm~\ref{accextra} only checks the decrease of $||X-W_nH_y||_F$ where $(W_n,H_y)$ is not necessarily feasible for $hp \geq 2$. The reason for computing $||X-W_nH_y||_F$ 
and not $||X-W_nH_n||_F$ is threefold: 
(i)~$W_n$ was updated according to $H_y$, 
(ii)~it gives the algorithm some degrees of freedom to possibly increase the objective function, 
(iii)~it is computationally cheaper because computing $||X-W_nH_n||_F$ would require $O(mnr)$ operations instead of $O(mr^2)$ (see Remark~\ref{comperr}). 
%The reason $H_y$ is reset as $H_n$ and not $H$ is because $H_n$ was optimized for $W_y$ which had previously decreased the error. However, this choice has a very minor impact on the performance of Algorithm~\ref{accextra}. 

In order to guarantee the objective function to decrease, a possible way is to require $H_y$ to be nonnegative by projecting it to the nonnegative orthant: this variant corresponds to $hp=3$.  
In that case, the solution $(W_n,H_y)$ is a feasible solution for which the objective function is guaranteed to decrease at least every second step. In fact, when the error increases, Algorithm~\ref{accextra} reinitializes the extrapolation sequence $(W_y,H_y)$ using $(W,H)$ (step~17 of Algorithm~\ref{accextra}) and the next step is a standard NNLS update. Therefore, since the objective function is bounded below, there exists a converging subsequence of the iterates. Proving convergence to a stationary points is an open problem, and an important direction for further research. 
We believe it would be particularly interesting to investigate the convergence of the extrapolation scheme applied on block coordinate descent in the non-convex case. 

To summarize, 
using the extrapolation of $H$ after the update of $W$ ($hp=1$) or using the projection of $H_y$ onto the feasible set ($hp=3$) is more conservative but guarantees the objective function to decrease (at least every second step). As we will see in the numerical experiments, these two variants perform in general better than with $hp=2$.

\begin{remark}[Computation of the error] \label{comperr}
To compute the error $||X-W_nH_y||_F^2$ in step~15 of Algorithm~\ref{accextra} (and in step~1), it is important to take advantage of previous computations and not compute $W_nH_y$ explicitly (which would be impractical for large sparse matrices). For simplicity, let us denote $W = W_n$ and $H = H_y$. 
We have 
\begin{align*}
||X-WH||_F^2 
& = \langle X , X \rangle - 2 \langle X, WH \rangle + \langle WH, WH \rangle \\
& = ||X||_F^2 - 2 \langle W, X H^T \rangle + \langle W^TW, HH^T \rangle . 
\end{align*} 
The term $||X||_F^2$ can be computed once, 
the term $\langle W, X H^T \rangle$ can be computed in $O(m r)$ operations since $M H^T$ is computed within the NNLS update of $W$, 
and the term $\langle W^TW, HH^T \rangle$ requires $O(m r^2)$ since $HH^T$ is also computed within the NNLS update of $W$. 
In fact, all algorithms for NNLS we know of need to compute $X H^T$ and $HH^T$ when solving for $W$, because the gradient of $||X-WH||_F^2$ with respect to $W$ is $2(WHH^T-XH^T)$. 
\end{remark}

\section{Choice of the extrapolation parameters $\beta_k$'s} \label{sec:beta}

In this section, we propose a strategy to choose the $\beta_k$'s.  
First, let us explain why it does not work well to use line search. Let us focus on the update of $W$ (a similar argument holds for $H$). 
We have 
\[
W_y = W_y(\beta) = W_n + \beta (W_n - W), 
\]
where $W_n$ is an approximate solution of $\min_{W \geq 0} ||X - WH_y||_F$ (in the case of ANLS, it is an optimal solution). %Note that for $\beta = -1$, we have $W_y = W$ and for $\beta = $, $W_y = W_n$. 
The optimal $\beta$ can be computed in close form as follows 
\[
\beta^* 
= \argmin_{\beta} ||X - W_y(\beta) H_y||_F^2 = 
\frac{ \langle X - W_n H_y ,(W_n - W) H_y \rangle }{ ||(W_n - W)H_y||_F^2 } .  
\]
We have observed that $\beta^*$ is close to zero for most steps of Algorithm~\ref{accextra} ($\beta^*$ is not always close to zero--even when using ANLS--because $W_y$ is not necessarily nonnegative), especially when the algorithm has performed several iterations and reached the neighbourhood of a stationary point. 
The reason is that $W_n$ was optimized to minimize the objective function. 
Hence, in the following, we propose another strategy to choose the $\beta_k$'s. It will increase the objective function in most cases (that is, $||X - W_y(\beta) H_y||_F^2 > ||X - W_y(0) H_y||_F^2$) but will allow a larger decrease of the objective function at the next step.  
Note that this is the reason why we check whether the error has decreased only after the update of $H$ because otherwise the acceleration would not be possible (only a small $\beta$ would be allowed in that case).

\paragraph{Strategy for updating the $\beta_k$'s} %\label{strat2} 

In this paper, since we are applying the extrapolation scheme to a non-convex problem using coordinate descent, there is, as far as we know, no a priori theoretically sound choice for the $\beta_k$'s. 
For this reason, we consider a very simple scheme described in Algorithm~\ref{betaupdate}. 

\algsetup{indent=2em}
\begin{algorithm}[ht!]
\caption{Update of the $\beta_k$'s \label{betaupdate}} 
\begin{algorithmic}[1]
\REQUIRE Parameters $1 < \bar{\gamma} < \gamma <  \eta$, $\beta_1 \in (0,1)$.
\STATE Set $\bar{\beta} = 1$.
\IF{ the error decreases at iteration $k$ } 
	\STATE Increase $\beta_{k+1}$: $\beta_{k+1} = \min( \bar{\beta}, \gamma \beta_{k} )$.
	\STATE Increase $\bar{\beta}$: $\bar{\beta} = \min\left(1, \bar{\gamma} \bar{\beta} \right)$. 
\ELSE {} 
	\STATE Decrease $\beta_{k+1}$: $\beta_{k+1} = \beta_{k+1} = \beta_{k}/\eta$.
	\STATE Set $\bar{\beta} = \beta_{k-1}$. 
\ENDIF  
\end{algorithmic} 
\end{algorithm} 

It works as follows. Let us assume there exists a hidden optimal value for the $\beta_k$'s, like in the strongly convex case where $\beta_k$ should ideally be equal to $\frac{1-\sqrt{\mu/L}}{1+\sqrt{\mu/L}}$~\cite{nesterov2013introductory, o2015adaptive}, where $\mu$ is the strong convexity parameter of the objective function and $L$ is the Lipschitz constant of its gradient. 
%Note that initially, the algorithm takes rather large steps since usually the initialization $(W,H)$ is far from being locally optimal, hence it makes sense for $\beta$ to be chosen smaller. 
It starts with an initial value of $\beta_0 \in [0, \bar{\beta}]$, and an upper bound $\bar{\beta}=1$.   
As long as the error decreases, it increases the value of $\beta_{k+1}$ by a factor $\gamma$ taking into account the upper bound, that is, $\beta_{k+1} = \min(\gamma \beta_k, \bar\beta)$. It also increases the upper bound by a factor $\bar{\gamma} < \gamma$ if it is smaller than one, that is, $\bar\beta = \min(\gamma \bar\beta,1)$. 
The usefulness of $\bar{\beta}$ is to keep in memory the last value of $\beta_k$ that allowed decrease of the objective function which is used as an upper bound for $\beta_k$. 
However, because the landscape of the objective function may change, $\bar{\beta}$ is slightly increased by a factor $\bar{\gamma} < \gamma$ at each step, as long as the error decreases. 
When the error increases, $\beta_{k+1}$ is reduced by a factor $\eta > \gamma$ and the upper bound $\bar\beta$ is set to the previous value of $\beta$ that allowed decrease, that is, $\beta_{k-1}$.

\begin{remark}
We have also tried to mimic the choice of the $\beta_k$'s from convex optimization~\cite{nesterov2013introductory} but it did in general perform worse than the simple choice presented here. 
\end{remark}

\section{Numerical Experiments} \label{numexp}  

In this section, we show the efficiency of the extrapolation scheme, that is, Algorithm~\ref{accextra},  
to accelerate the NMF algorithms ANLS and A-HALS. All tests are preformed using Matlab
R2015a on a laptop Intel CORE i7-7500U CPU @2.9GHz 24GB RAM. 
The code is available from \url{https://sites.google.com/site/nicolasgillis/code}. 

\subsection{Data sets} 

We will use the same data sets as in~\cite{gillis2012accelerated} as they are among the most widely used ones in the NMF literature; see Tables~\ref{imdat} and \ref{dtm}. 
The image data sets represent facial images and are dense matrices. 
The document data sets are sparse matrices. 
\begin{center}
\begin{table}[h!]
\begin{center}
\caption{Image datasets. }
\label{imdat}  
\begin{tabular}{|c|c|c|c|c|}
\hline 
Name &             $\#$ pixels &  $m$  &  $n$ &  $r$   \\ \hline \hline
ORL$^1$    &   $112 \times 92$  & 10304 & 400 & 40    \\ 
Umist$^2$  &   $112 \times 92$  & 10304 & 575 & 40   \\ 
CBCL$^3$ &  $19 \times 19$  & 361  & 2429  & 40 \\ 
Frey$^2$  &   $28 \times 20$ & 560 & 1965  & 40  \\ 
\hline
\end{tabular} \\
\begin{flushleft}

\footnotesize
$^1$ \url{http://www.cl.cam.ac.uk/research/dtg/attarchive/facedatabase.html}\\
$^2$ \url{http://www.cs.toronto.edu/~roweis/data.html}\\
$^3$ \url{http://cbcl.mit.edu/software-datasets/FaceData2.html} 
\end{flushleft}
\end{center}
\end{table}
\end{center}
\begin{center}
\begin{table}[h!] 
\begin{center}
\caption{Text mining datasets \cite{ZG05} (sparsity is given in $\%$: $100*\#\text{zeros}/(mn)$).} 
\label{dtm}
\begin{tabular}{|c|c|c|c|c|c|} 
\hline
Name &   $m$  &  $n$ &  $r$ &  $\#\text{nonzero}$ & sparsity   \\ \hline \hline
classic &  7094   & 41681 &  20 & 223839  & 99.92 \\  
sports &    8580 & 14870 &   20 & 1091723 & 99.14 \\ 
reviews &   4069  & 18483 &  20 & 758635  & 98.99\\ 
hitech &    2301 & 10080 &   20 & 331373 & 98.57\\ 
ohscal &   11162  & 11465 &  20 & 674365 & 99.47  \\ 
la1 &    3204 & 31472 &   20 & 484024  & 99.52 \\ 
\hline
\end{tabular}
\end{center}
\end{table}
\end{center}

We also consider two types of synthetic data sets. 
For the first one, which we will refer to as the low-rank synthetic data sets, we generate each entry of $W$ and $H$ using the uniform distribution in $[0,1]$ and compute $X = WH$. 
For each experiment, we will generate 10 such matrices and report the average results. For the second one, which we will refer to as the full-rank synthetic data sets, we simply generate each entry of $X$ uniformly at random in [0,1] so that $X$ is a full rank matrix. In both cases, we will use $m=n=200$ and $r=20$.

\subsection{Experimental set up} 

In all cases, we will report the average error over 10 random initializations, where the entries of the initial matrices $W$ and $H$ are chosen uniformly at random in the interval $[0,1]$. 
To compare the solutions generated by the different algorithms, we follow the strategy from~\cite{gillis2012accelerated}, that is, we report the relative error to which we subtract the lowest relative error obtained by any algorithm with any initialization (denoted $e_{\min}$). 
Mathematically, given the solution $(W^{(k)},H^{(k)})$ obtained at the $k$th iteration, we will report 
\begin{equation} \label{relerr}
E(k) = \frac{||X - W^{(k)} H^{(k)}||_F}{||X||_F} - e_{\min}. 
\end{equation}
%where $e_{\min}$ is the smallest relative error obtained among all algorithms and all initializations. 
For the low-rank synthetic data sets, we use $e_{\min} = 0$. 

Using $E(k)$ instead of $||X - W^{(k)} H^{(k)}||_F$ has several advantages: 
(i)~it allows to take meaningfully the average results over several data sets, 
and  
(ii)~it provides a better visualization both in terms of initial convergence and in terms of the quality of the final solutions computed by the different algorithms. The reason is that $E(k)$ converges to zero for the algorithm that was able to compute the best solution which allows us to use a logarithmic scale. 
%However, since we will take averages (across different data sets and initializations), the curves on the plots will not necessarily converge to the 0. 
%Also, when comparing NMF algorithms that may converge to different local minima, 
%comparing the values of $E(k)$ especially make sense when looking at the speed of convergence. 

\subsection{Tuning parameters: preliminary numerical experiments} 

Before we compare the two NMF algorithms (ANLS and A-HALS) and their extrapolated variants, 
we run some preliminary numerical experiments in order to choose reasonable values for 
the parameter of Algorithm~\ref{accextra} ($hp$) and the parameters to update $\beta_k$. 

As we will see, the extrapolation scheme performs rather differently for ANLS (that computes an optimal solution of the subproblems) and A-HALS (that computes an approximate solution using a few steps of coordinate descent). 
It will also perform rather differently depending on the value of $hp$, 
while it will less sensitive to the values of $\beta_0$, $\gamma$, $\bar\gamma$ and $\eta$ as long as these values are chosen in a reasonable range. 
%In fact, for ANLS, it will be more effective to perform the extrapolation only on $W$ (that is, use $hp = 1$) while for A-HALS, it will be more effective to perform it on both $W$ and $H$ (that is, $hp=2,3)$. 

In the next section, we will run the different variants with the following parameters: 
 $\beta_0 = 0.25, 0.5, 0.75$, 
 $\eta = 1.5, 2, 3$, 
 $(\gamma, \bar{\gamma}) = (1.01,1.005), (1.05,1.01), (1.1,1.05)$. 
For each experiment, we will not be able to display the curve for each extrapolated variants (there would be too many, 82 in total: $3^4$ and the original algorithm). Therefore, for each value of $hp$, we will only display the variant corresponding to the parameters that obtained the smallest final average error (best) and the largest final average error (worst). This will be interesting to observe the sensitivity of Algorithm~\ref{accextra} to the way $\beta_k$ is updated.

%ANLS low-rank: hp=1, $\beta_0 = 0.5$, $\eta = 1.5$, $\gamma = 1.1$. 
%ANLS full rank: hp=1, $\beta_0 = 0.25$, $\eta = 2$, $\gamma = 1.1$. 
%A-HALS low-rank: hp=2, $\beta_0 = 0.75$, $\eta = 1.5$, $\gamma = 1.01$. 
%(Second: hp=3, $\beta_0 = 0.5$, $\eta = 1.5$, $\gamma = 1.01$.) 
%A-HALS full rank: hp=1, $\beta_0 = 0.75$, $\eta = 2$, $\gamma = 1.05$. (hp = 2,3 also works well.) \\

\subsubsection{Extrapolated ANLS (E-ANLS)} 

\begin{figure}[ht!]
\begin{center}
\begin{tabular}{cc}
\includegraphics[width=0.45\textwidth]{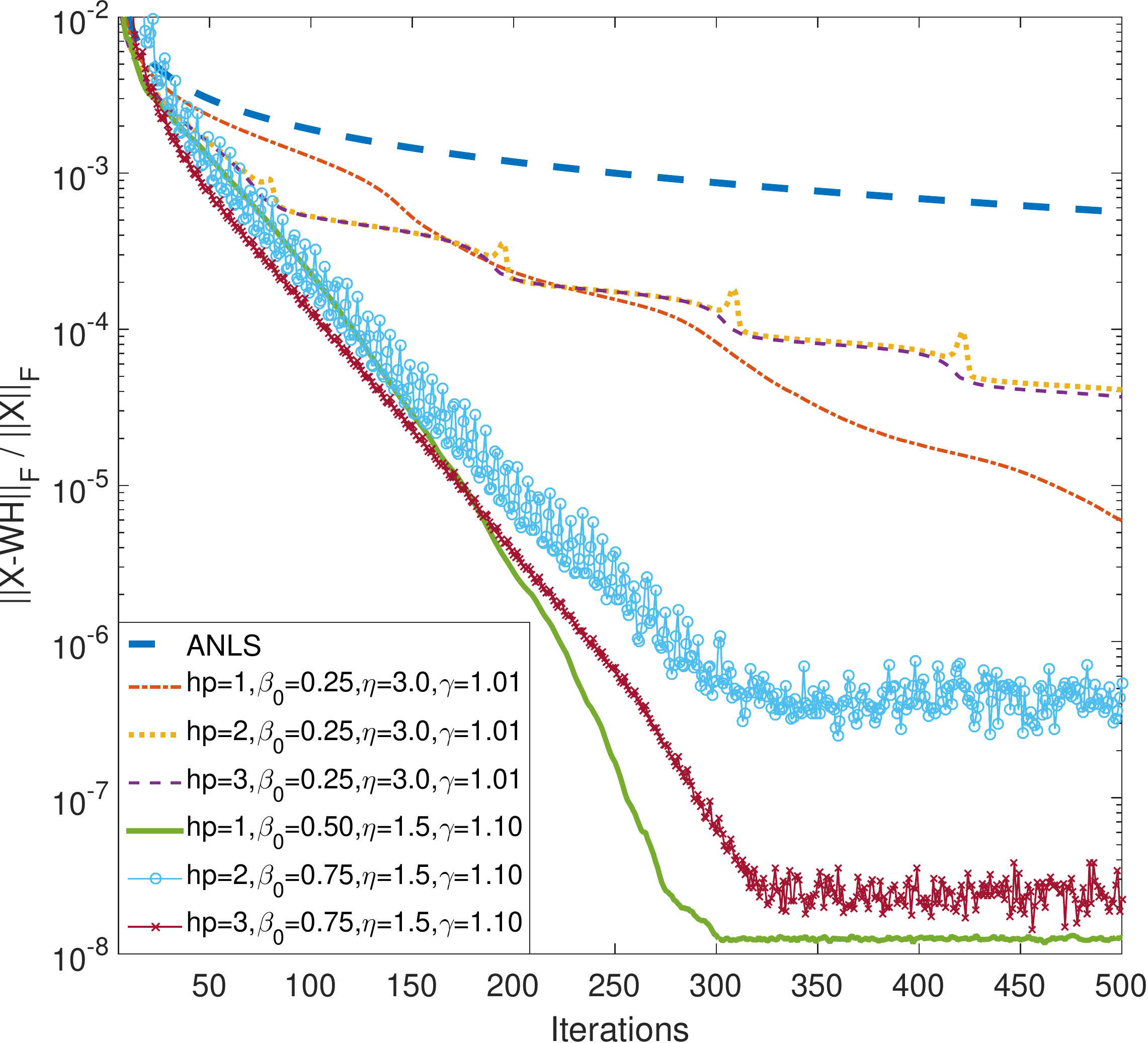}  & 
\includegraphics[width=0.45\textwidth]{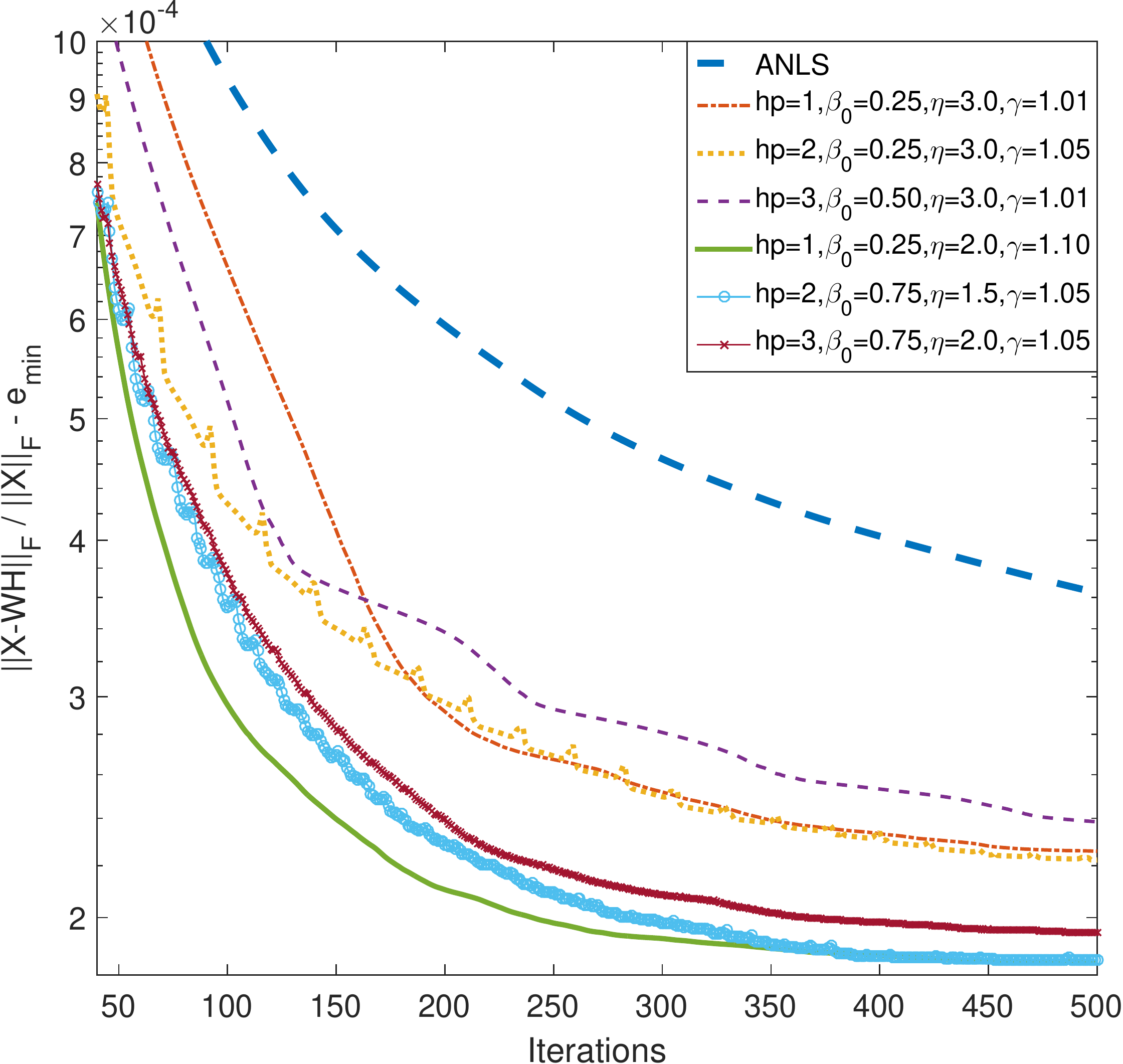} \\ 
\includegraphics[width=0.45\textwidth]{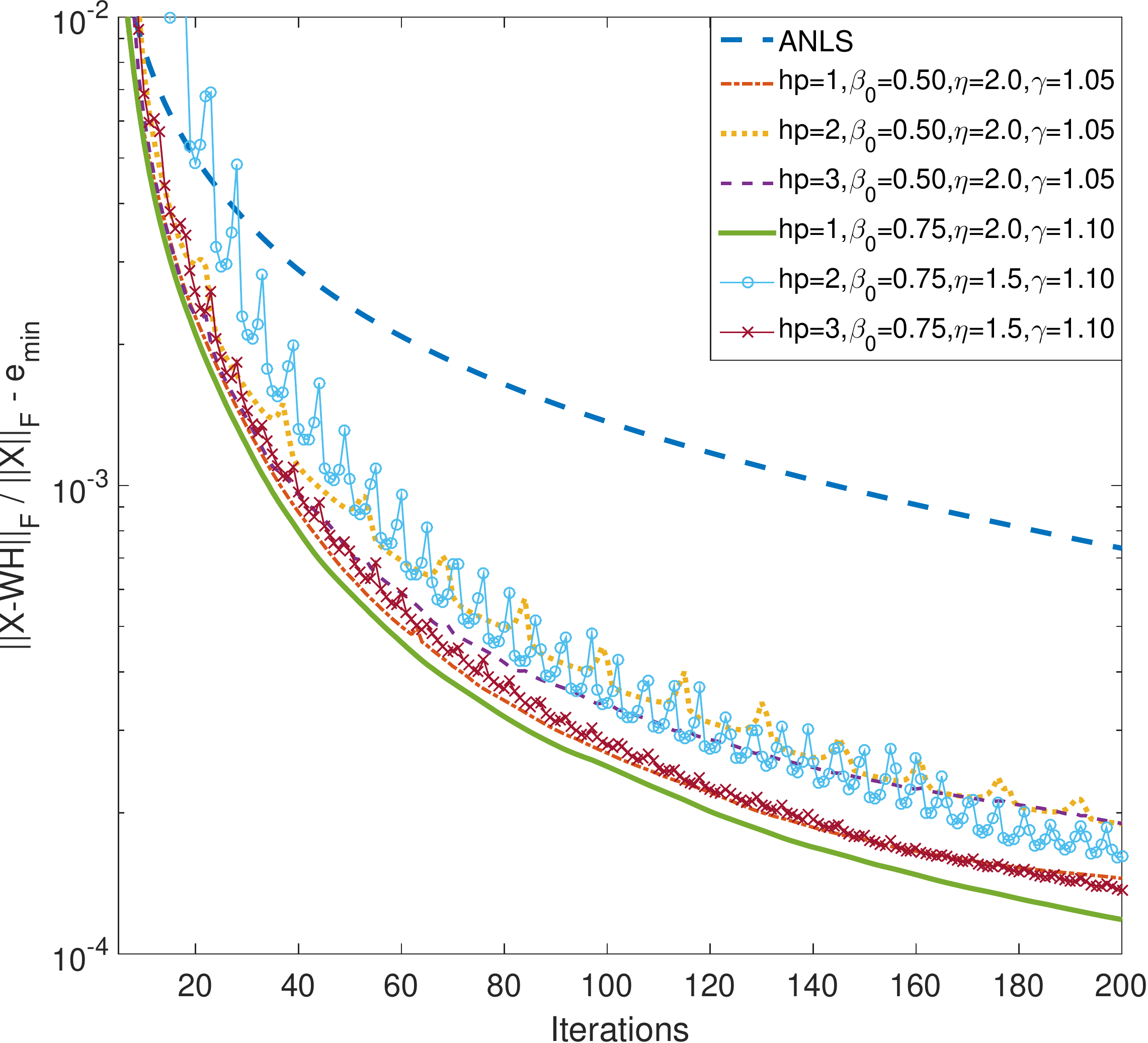}  & 
\includegraphics[width=0.45\textwidth]{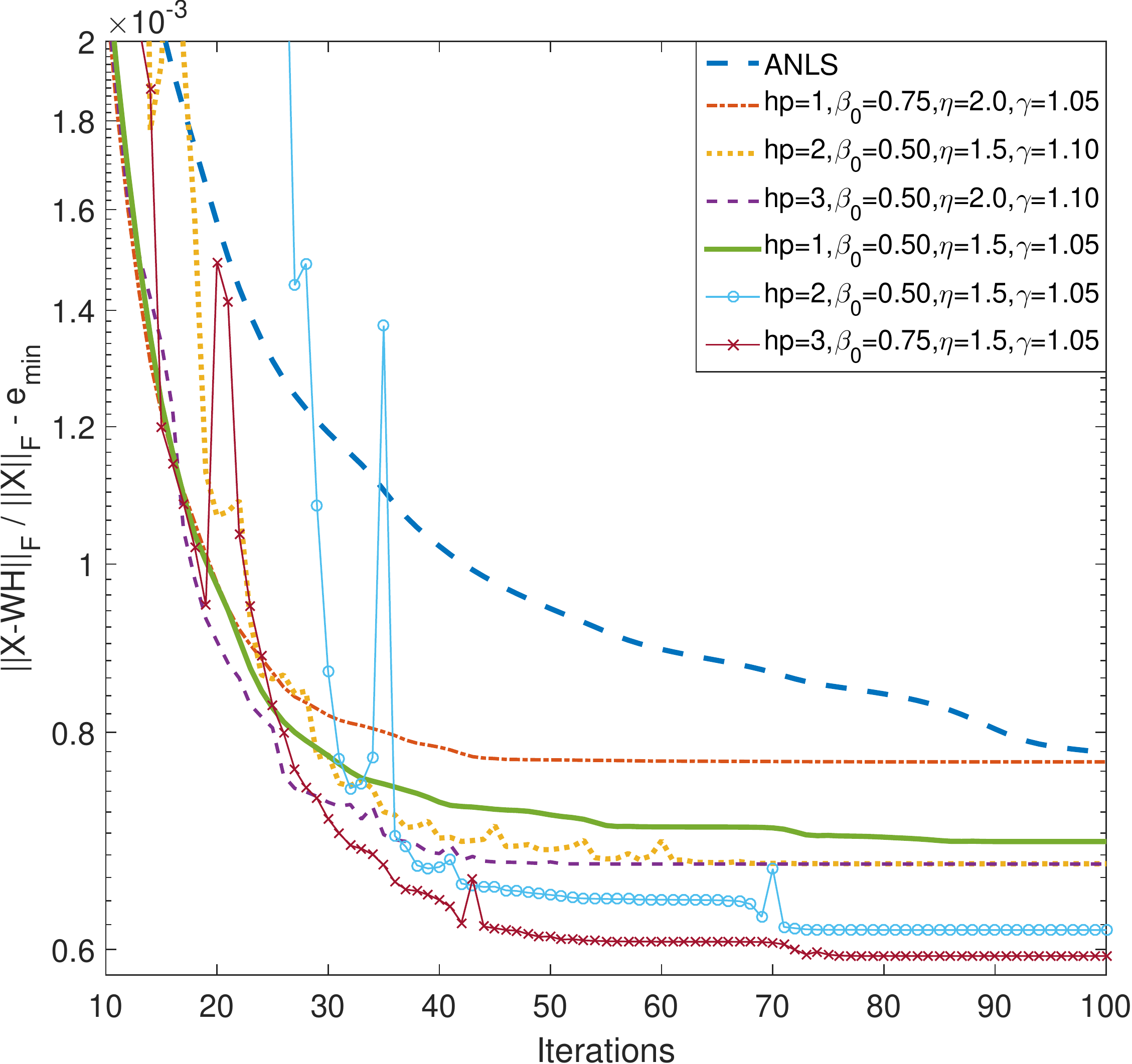} \\ 
\end{tabular}
\caption{
Extrapolation scheme applied to ANLS on the low-rank (top left) and full-rank (top right) synthetic data sets, and the image (bottom left) and document (bottom right) real data sets.  
For each value of $hp$, we display the corresponding best and worst performing variant.  
The curves are the average value of~\eqref{relerr} among the different data sets and initializations.  \label{fig:synt_anls}}
\end{center}
\end{figure}

The top two plots of Figure~\ref{fig:synt_anls} show the evolution of the average of the error measure defined in~\eqref{relerr} for the low-rank and full-rank synthetic data sets. 
We observe the following:
\begin{itemize}
\item For all the values of the parameters, E-ANLS outperforms ANLS. 

\item For the low-rank synthetic data sets,  
E-ANLS with $hp=1$ and well-chosen parameters for the update of $\beta_k$ (e.g., $\beta_0=0.5,\eta=1.5,\gamma=1.1$) performs extremely well and is able to identify solutions with very small relative error 
($\approx 10^{-8}$ in average). 
In fact, the original ANLS algorithm would not be able to compute such solutions even within several thousand iterations. 

\item For the full-rank synthetic data sets, E-ANLS variants with $hp$ equal to 1, 2 or 3 perform similarly although choosing $hp=1$ allows a slightly faster initial convergence. 

\item The best value for $\gamma$ is either 1.05 or 1.10. 
The best value for $\eta$ is either 1.5 or 2 (3 being always the worst). 
The algorithm is not too sensitive to the initial $\beta$ as it is quickly modified within the iterations but $\beta_0 = 0.25$ clearly provides the worst performance in most cases.   

\end{itemize}

We now perform the same experiment on image and document data sets except with fewer parameters 
(we do not use $\gamma = 1.01$, $\eta = 3$, $\beta_0 = 0.25$) in order to reduce the computational load. 
The bottom two plots of Figure~\ref{fig:synt_anls}  show the evolution of the average of the error measure defined in~\eqref{relerr} for the image and document data sets. 
% \begin{figure}[ht!]
%\begin{center}
%\begin{tabular}{cc}
%\includegraphics[width=0.45\textwidth]{figures/anls_images}  & 
%\includegraphics[width=0.45\textwidth]{figures/anls_text} \\ 
%\end{tabular}
%\caption{
%Extrapolation scheme applied with ANLS on the image (left) and document (right) data sets. 
%For each value of $hp$, we display the corresponding best and worst performing variant.   \label{fig:real_anls}}
%\end{center}
%\end{figure}

We observe the following:
\begin{itemize}
\item As for synthetic data sets, E-ANLS outperforms ANLS for all the values of the parameters. 

\item Since we have removed the values of the parameters performing worst, the gap between the best and worst  variants of E-ANLS is reduced. 

\item For the image data sets, the variant with $hp=1$ performs best although the variants with $hp=2,3$ do not perform much worse. 

 \item For the document data sets, 
 the variants with $hp=2,3$ perform best (in terms of final error). 
 This can be explained by the fact that NMF problems for sparse matrices are more difficult as there are more local minima with rather different objective function values (see Section~\ref{sec:text} for more numerical experiments). Hence the final error reports the algorithm that found the best solution in most of the 60 cases (6 data sets, 10 initializations per data set).   
In terms of speed of convergence, most E-ANLS variants behave similarly (converging within 80 iterations, while ANLS has not converged within 100 iterations).   

\end{itemize}

In the final numerical experiments, 
we will use $\beta_0 = 0.5$, $\eta = 1.5$ and $(\gamma,\bar{\gamma}) = (1.1,1.05)$ for E-ANLS. 
We will keep both variants $hp=1,3$.

\subsubsection{Extrapolated A-HALS (E-A-HALS)} 

\begin{figure}[ht!]
\begin{center}
\begin{tabular}{cc}
\includegraphics[width=0.45\textwidth]{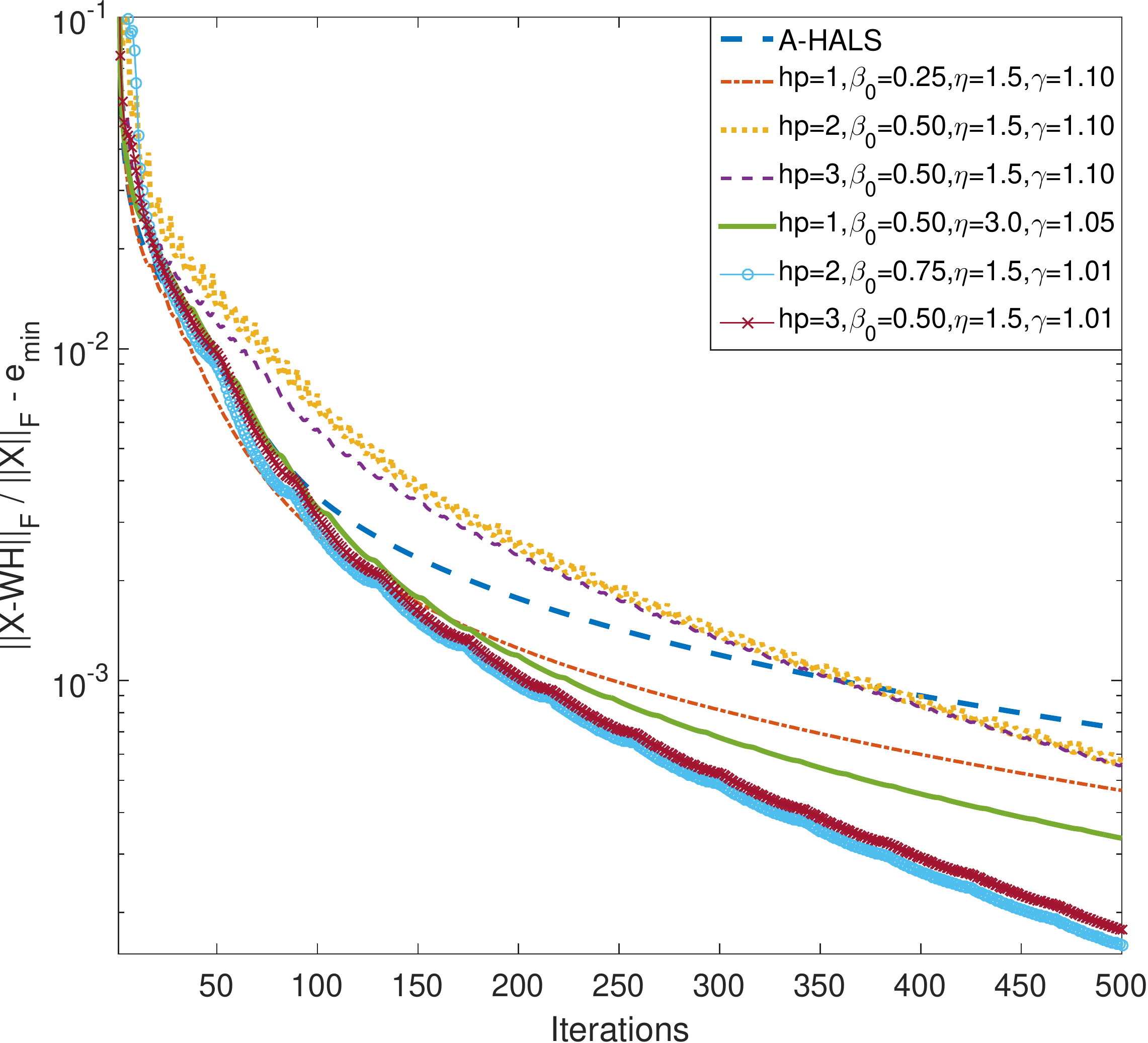}  & 
\includegraphics[width=0.45\textwidth]{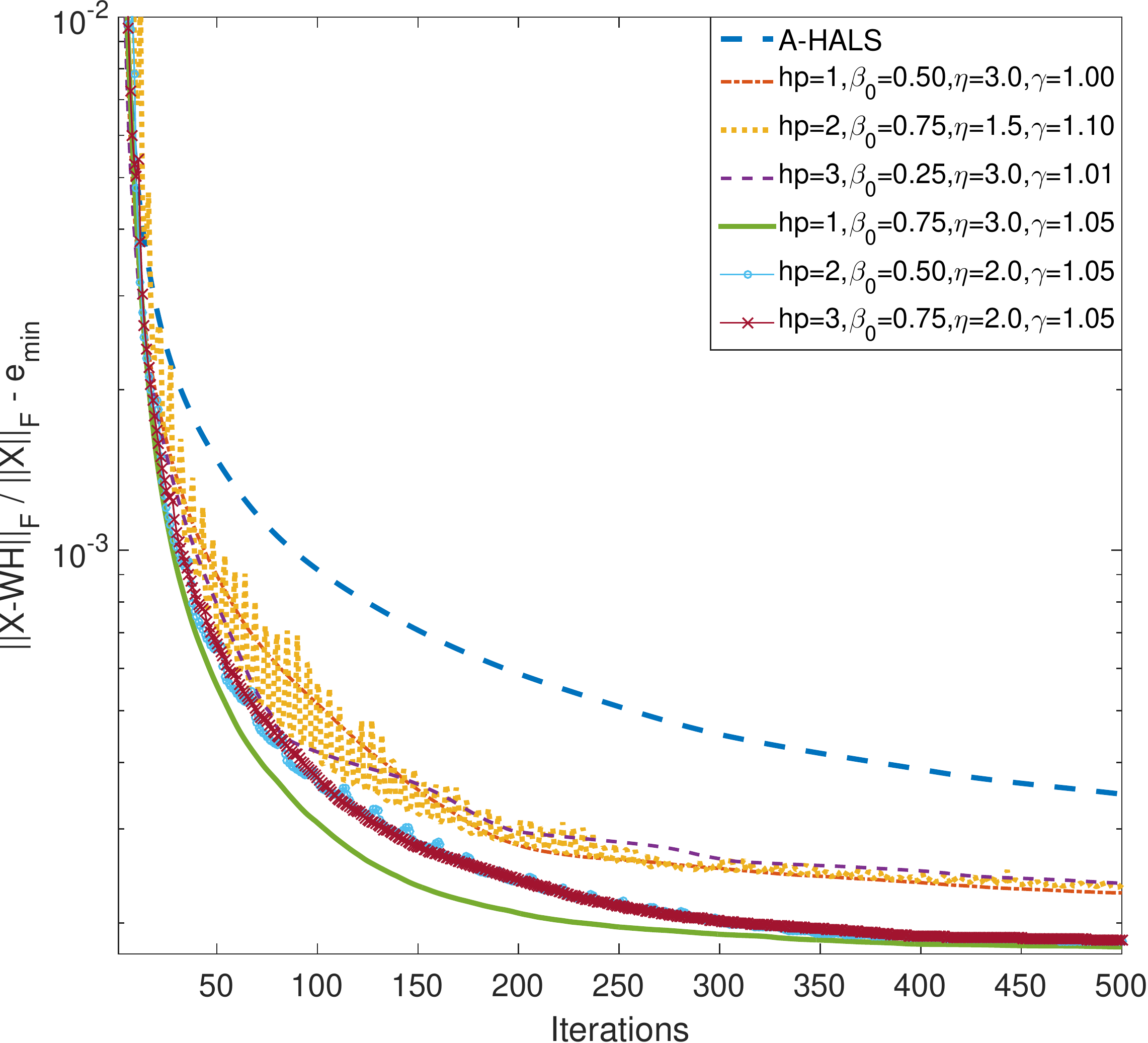} \\ 
\includegraphics[width=0.45\textwidth]{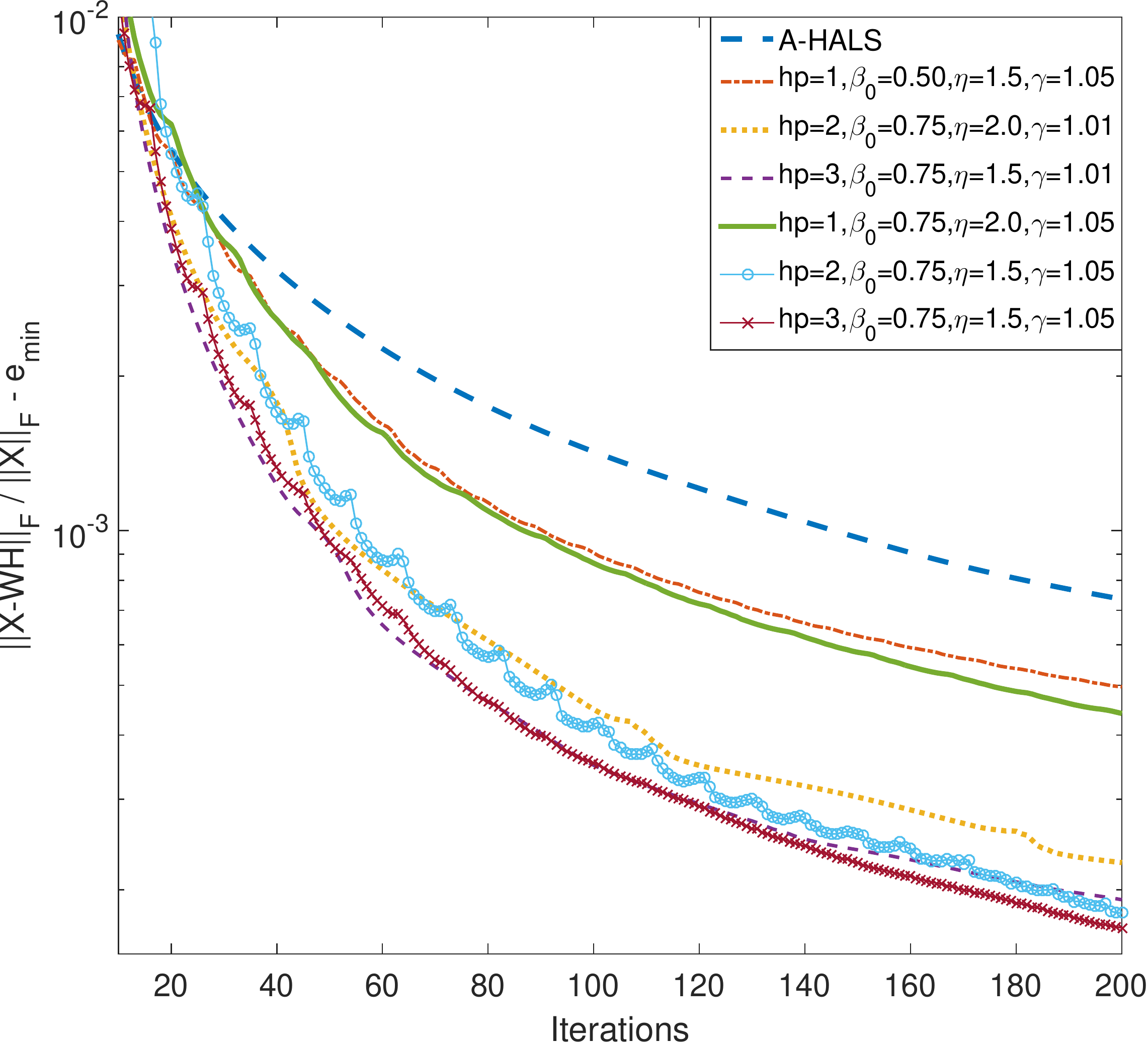}  & 
\includegraphics[width=0.45\textwidth]{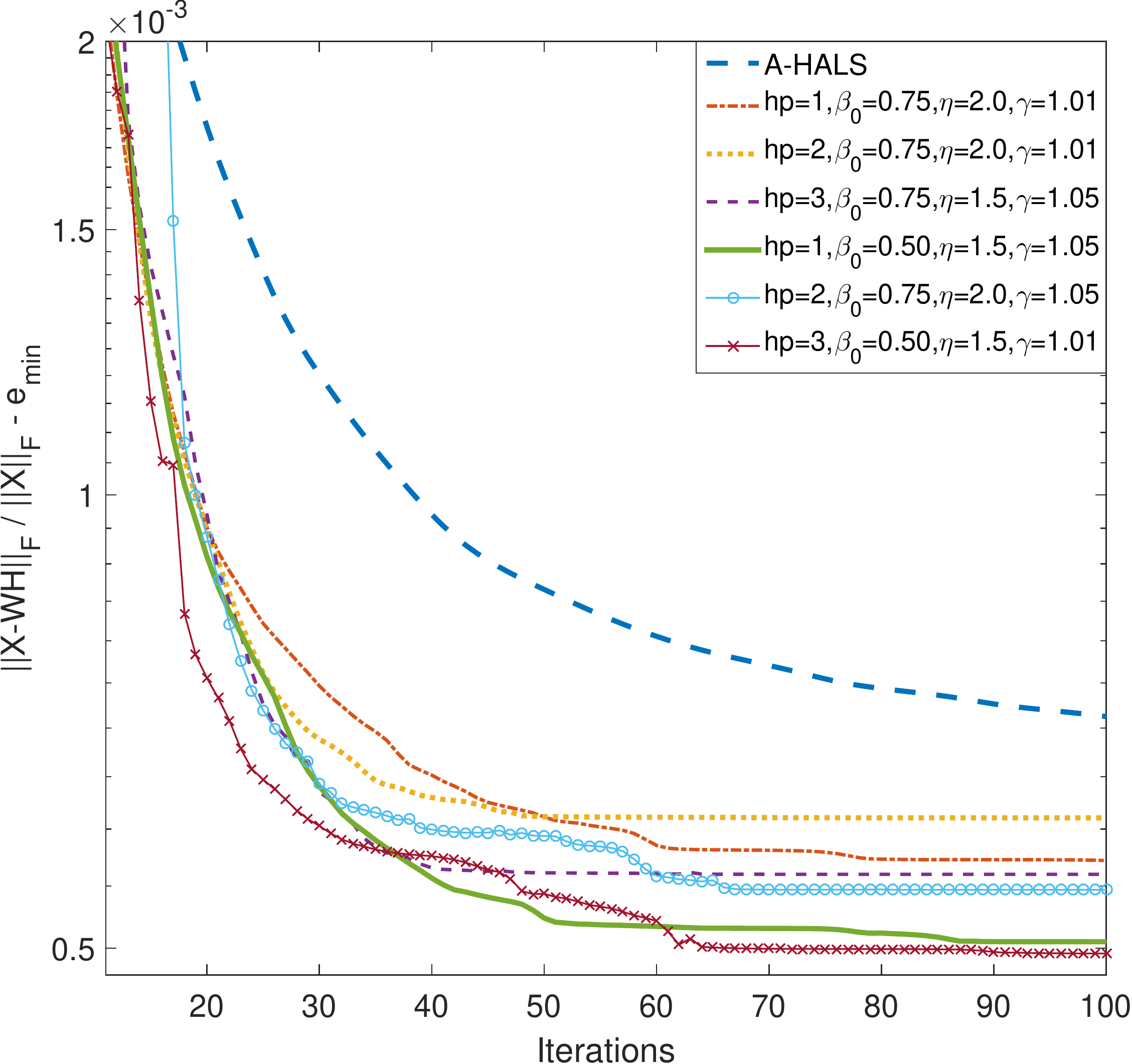} \\ 
\end{tabular}
\caption{
Extrapolation scheme applied with A-HALS on the low-rank (top left) and full-rank (top right) synthetic data sets, and on the image (bottom left) and document (bottom right) data sets.  
For each value of $hp$, we display the corresponding best and worst performing variant.   \label{fig:synt_hals}}
\end{center}
\end{figure}

The top two plots of Figure~\ref{fig:synt_hals} show the evolution of the average of the error measure defined in~\eqref{relerr} for the low-rank and full-rank synthetic data sets. For these experiments, we have also tested the value 
$(\gamma,\bar{\gamma})=(1.005,1.001)$ (as we will see that smaller value of these parameters perform better). 
We observe the following:
\begin{itemize}

\item For the low-rank synthetic data,  with $hp=2,3$ and well-chosen parameters for the update of $\beta$ 
(e.g., $\beta_0=0.50,\eta=1.5,\gamma=1.01$), E-A-HALS performs much better than A-HALS. (Note however that it is not able to find solutions with error as small as E-ANLS within 500 iterations.) 

\item For the full-rank synthetic data, the variant with $hp=1$ performs slightly better although the final solutions of the three extrapolated variants have similar error. 

\item The best value for $\gamma$ is either 1.01 or 1.05, smaller than for E-ANLS. This can be explained by the fact that A-HALS does not solve the NNLS subproblems exactly and the extrapolation should not be as aggressive as for ANLS. As for E-ANLS, E-HALS is not too sensitive to the parameters $\eta$ and $\beta_0$. 

\end{itemize}

We now perform the same experiment on image and document data sets except with fewer parameters 
(we do not test $\gamma = 1.005, 1.1$, $\eta = 3$, $\beta_0=0.25$).  
The bottom two plots of Figure~\ref{fig:synt_hals}  show the evolution of the error measure defined in~\eqref{relerr} for the image and document data sets.  
% \begin{figure}[ht!]
%\begin{center}
%\begin{tabular}{cc}
%\includegraphics[width=0.45\textwidth]{figures/hals_images}  & 
%\includegraphics[width=0.45\textwidth]{figures/hals_text} \\ 
%\end{tabular}
%\caption{
%Extrapolation scheme applied with ANLS on the image (left) and document (right) data sets. 
%For each value of $hp$, we display the corresponding best and worst performing variant.   \label{fig:real_hals}}
%\end{center}
%\end{figure} 
We observe the following:
\begin{itemize}
%\item In both cases, the variant $hp=3$ performs slightly better than. 

\item For the image data sets, the variant $hp=1$ performs worse than $hp=2,3$ which perform similarly (in terms of speed of convergence).  

\item For the document data sets, we observe a similar behavior as for ANLS: all extrapolated variants converge much faster than HALS, but converge to different solutions (being in average less than 0.1\% away from the lowest relative error). 

\end{itemize}

In the final numerical experiments, 
we will use $\beta_0 = 0.5$, $\eta = 1.5$ and $(\gamma,\bar{\gamma}) = (1.01,1.005)$ for E-A-HALS.  
We will keep both variants $hp=1,3$.

\subsection{Extensive numerical experiments and comparison of E-ANLS and E-HALS} 

We now compare ANLS, A-HALS and their extrapolated variants on the same data sets. 
We also compare these algoirthms with the extrapolated alternating projected gradient method for NMF proposed by Xu and Yin~\cite{xu2013block} and referred to as PGM-MF.

\subsubsection{ Synthetic data sets }

\begin{figure}[ht!] 
\begin{center}
\begin{tabular}{cc}
\includegraphics[width=0.48\textwidth]{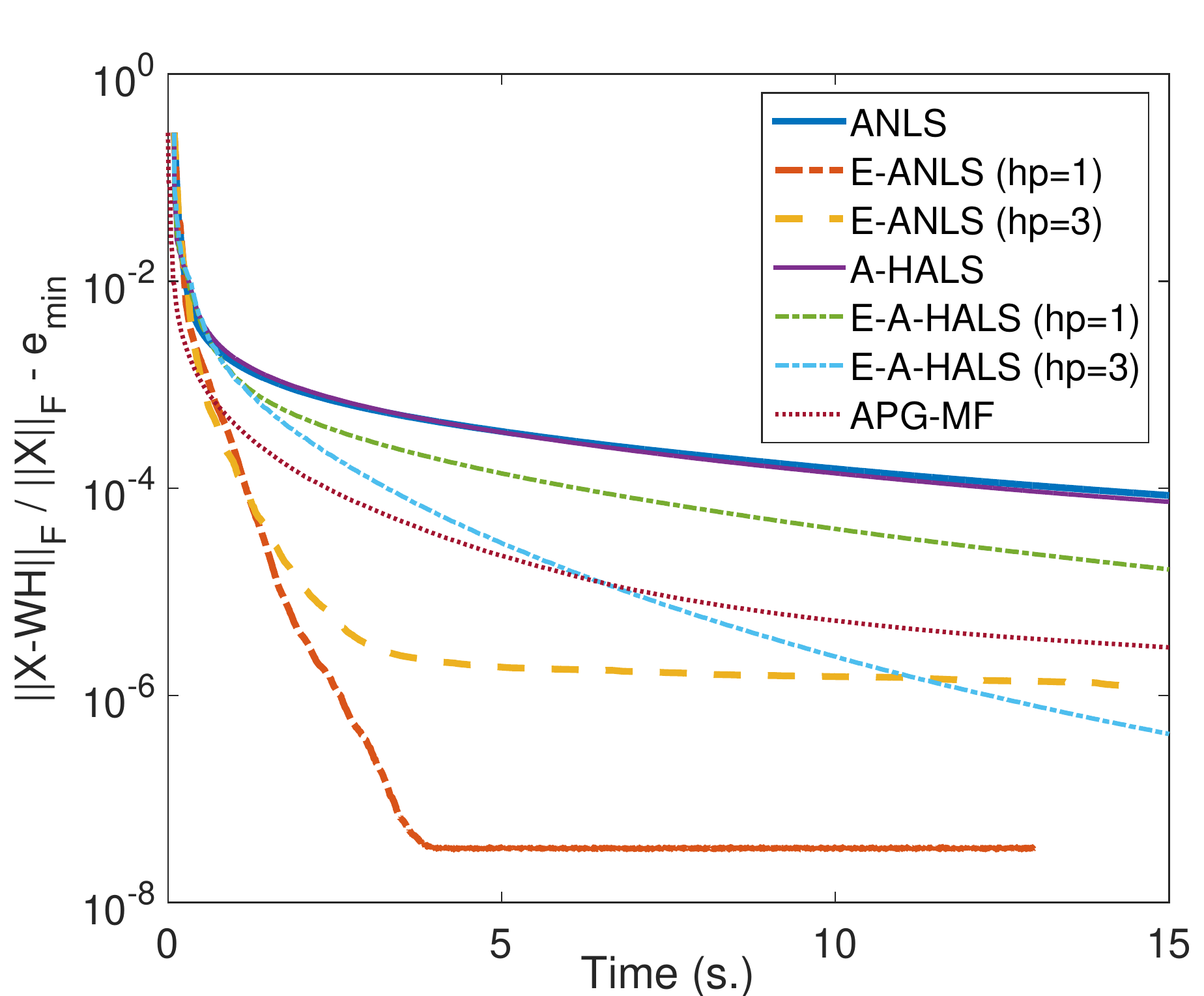}  & 
\includegraphics[width=0.48\textwidth]{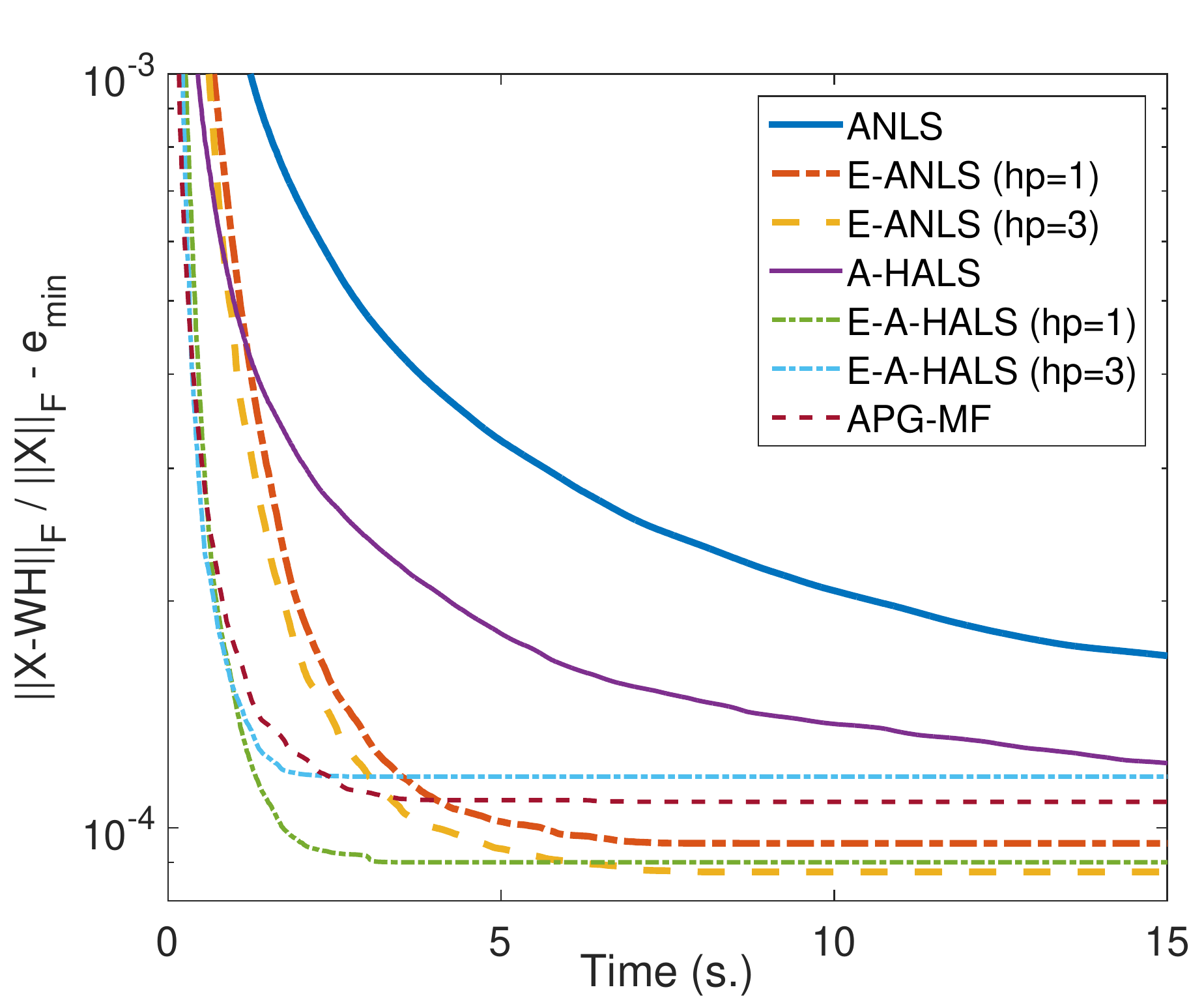} \\ 
\end{tabular}
\caption{
Average value of the error measure~\eqref{relerr} of ANLS, A-HALS and their extrapolated variants applied on low-rank (left) and full-rank (right) synthetic data sets.   \label{fig:randfinal}} 
\end{center}
\end{figure}

Figure~\ref{fig:randfinal} displays the evolution of the average error for the low-rank and full rank synthetic data sets, where the NMF algorithms were run for 20 seconds.  
Tables~\ref{tab:randlowrank} (resp.\@ \ref{tab:rand}) reports the average error, standard deviation and a ranking among the final solutions obtained by the different algorithms for the low-rank (resp.\@ full-rank) synthetic data sets. 
\begin{center}  
 \begin{table}[h!] 
 \begin{center} 
\caption{Comparison of the final relative error obtained by the NMF algorithms on the low-rank synthetic data sets: 
Average error, standard deviation and rankings among the 100 runs (100 data sets). 
The $i$th entry of the vector indicates the number of times the algorithm generated the $i$th best solution. (Observe that all algorithms are able to compute the best solution at least a few times: this happens when they compute an exact solution with $X=WH$.) 
\label{tab:randlowrank}} 
  \begin{tabular}{|c|c|c|c|c|c|c|} 
 \hline Algorithm &  mean $\pm$ std & ranking  \\ 
 \hline 
 ANLS &  $5.612\, 10^{-5} \pm 7.414\, 10^{-5}$ &  ( 0,  0,  0,  3,  7, 40, 50)   \\ 
E-ANLS ($hp=1$) &  $\mathbf{2.618\, 10^{-8} \pm 3.657\, 10^{-8}}$ &  (96,  4,  0,  0,  0,  0,  0)   \\ 
E-ANLS ($hp=3$) &  $1.207\, 10^{-6} \pm 1.162\, 10^{-5}$ &  (67, 24,  7,  1,  1,  0,  0)   \\ 
 A-HALS &  $4.547\, 10^{-5} \pm 6.299\, 10^{-5}$ &  ( 1,  0,  0,  4, 10, 41, 44)   \\ 
E-A-HALS ($hp=1$) &  $7.825\, 10^{-6} \pm 1.531\, 10^{-5}$ &  ( 3,  0,  6, 31, 41, 13,  6)   \\ 
E-A-HALS ($hp=3$) &  $1.181\, 10^{-7} \pm 3.679\, 10^{-7}$ &  (48,  8, 37,  7,  0,  0,  0)   \\ 
APG-MF &  $2.032\, 10^{-6} \pm 5.770\, 10^{-6}$ &  ( 0,  0,  3, 50, 41,  6,  0)   \\ 
\hline \end{tabular} 
 \end{center} 
 \end{table} 
 \end{center} 
 
 \begin{center}  
 \begin{table}[h!] 
 \begin{center} 
\caption{Comparison of the final relative error obtained by the NMF algorithms on the full-rank synthetic data sets: 
Average error, standard deviation and rankings among the 100 runs (10 data sets, 10 initializations each). 
The $i$th entry of the vector indicates the number of times the algorithm generated the $i$th best solution. 
\label{tab:rand}} 
 \begin{tabular}{|c|c|c|c|c|c|c|} 
 \hline Algorithm &  mean $\pm$ std & ranking  \\ 
 \hline 
 ANLS &              $0.423858 \pm 1.183\, 10^{-3}$ &  ( 4,  9,  9, 12, 22, 21, 23)   \\ 
E-ANLS ($hp=1$) &    $0.423795 \pm 1.161\, 10^{-3}$ &  (16, 18, 18,  9, 15, 10, 14)   \\ 
E-ANLS ($hp=3$) &    $\mathbf{0.423787 \pm 1.158\, 10^{-3}}$ &  (18, 11, 17, 21, 16,  9,  8)   \\ 
 A-HALS &            $0.423815 \pm 1.171\, 10^{-3}$ &  (18, 12, 11, 18, 13, 13, 15)   \\ 
E-A-HALS ($hp=1$) &  $0.423790 \pm 1.162\, 10^{-3}$ &  (17, 17, 16, 17, 15, 10,  8)   \\ 
E-A-HALS ($hp=3$) &  $0.423817 \pm 1.184\, 10^{-3}$ &  (12, 14, 16, 11, 11, 22, 14)   \\ 
APG-MF &             $0.423808 \pm 1.183\, 10^{-3}$ &  (16, 18, 13, 12,  8, 15, 18)   \\ 
\hline \end{tabular} 
 \end{center} 
 \end{table} 
 \end{center} 
 
 For low-rank synthetic data sets, these results confirm what we have observed previously: 
E-ANLS ($hp=1$) is able to significantly accelerate ANLS and obtain solutions with very small error extremely fast (in less than 4 seconds).  
The acceleration of HALS is not as important but it is significant.  
E-ANLS ($hp=1$) is able to obtain the best solutions in 96 out of the 100 runs while being always among the two best, while ANLS and A-HALS are among the worst ones in most cases. APG-MF never generates the best nor the second best solution. 
 
  For full-rank synthetic data sets, we observe that all algorithms obtain a similar final relative error (see Table~\ref{tab:rand}), all of them being in average around $0.01\%$ away from the best solution, and there is no clear winner between the extrapolated variants. In fact, there is a priori no reason to believe that an algorithm will converge to a better solution in general as NMF is a difficult non-convex optimization problem~\cite{vavasis2009complexity}. 
  %(For the low-rank synthetic data sets, the algorithms are able to recover an almost exact solution because these types of randomly generated NMF problems are not as hard.)  
  %--only ANLS and A-HALS were no able to converge within 20 seconds. 
  %It is important to recall the meaning of the error measure~\eqref{relerr}: it gives the difference of the relative error from $e_{\min}$ (the lowest error among the solutions found by all algorithms with all initializations) hence having a value around $2 \, 10^{-4}$ (like ) or $1.66 \, 10^{-4}$ (like ) is similar  
  %for example, E-A-HALS ($hp=3$) has the smallest average value, namely, , which means that in average over the 100 runs, 
  %it had relative error away  of $0.0166\%$. Hence, what is particularly important on to focus on is the initial speed of convergence. 
  In terms of speed of convergence, 
  E-A-HALS variants converge the fastest (about 3 seconds),  
   followed by APG-MF (about 4 seconds) and the E-ANLS variants (about 8 seconds), 
   while A-HALS and ANLS require more than 20 seconds.

\subsubsection{ Dense image data sets }  \label{sec:ima}

\begin{figure}[ht!]
\begin{center}
\begin{tabular}{cc}
\includegraphics[width=0.48\textwidth]{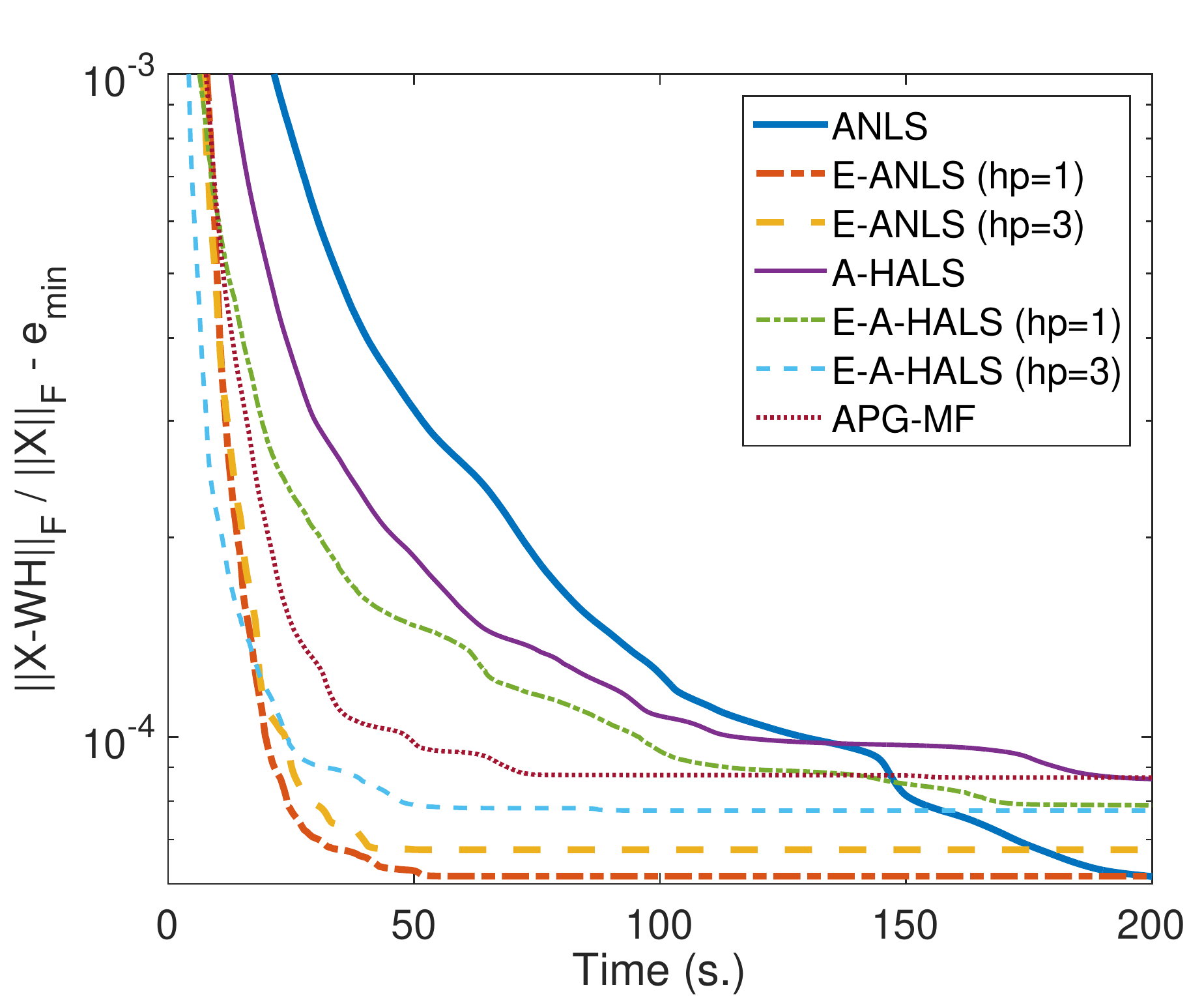}  & 
\includegraphics[width=0.48\textwidth]{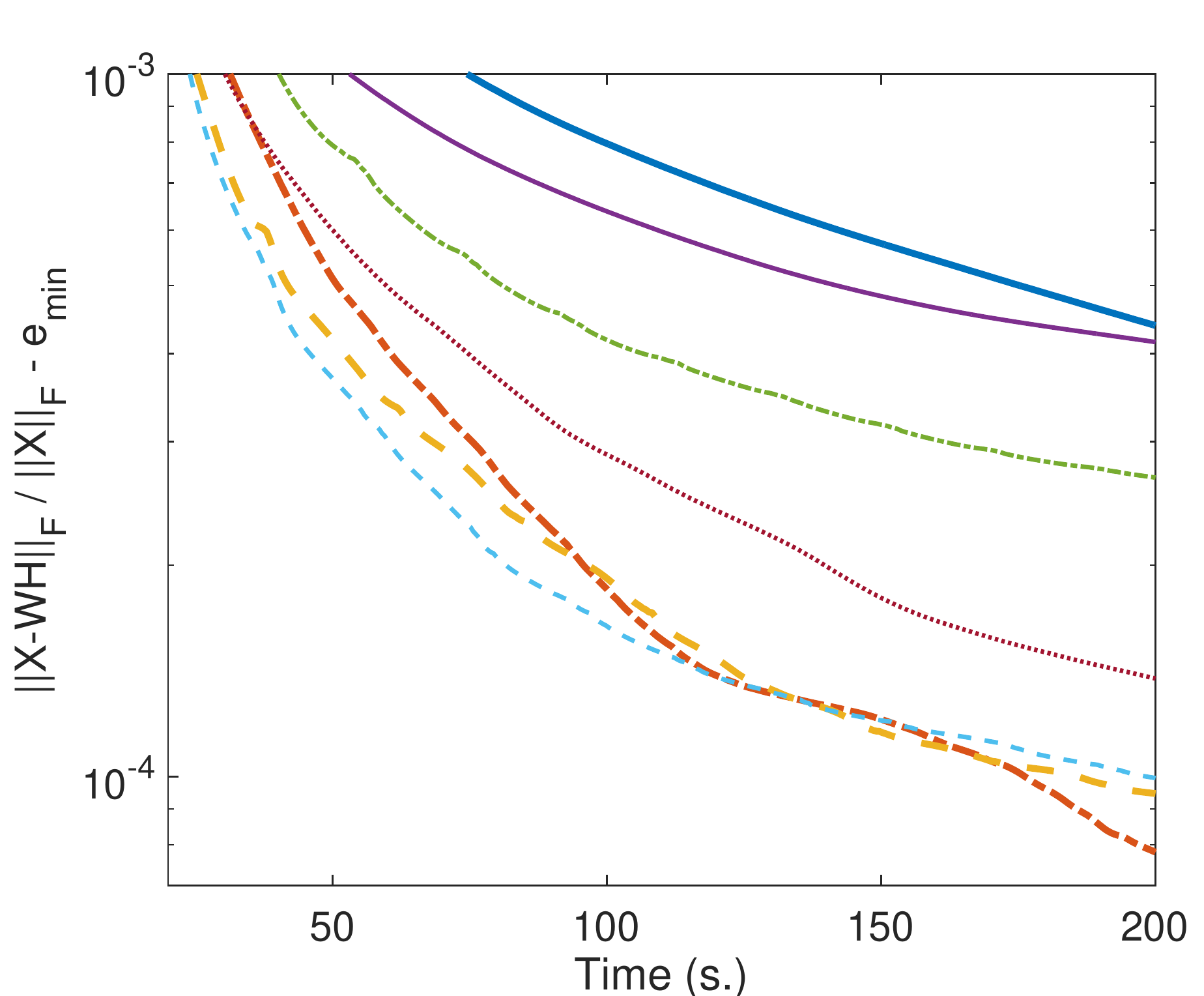} \\ 
\includegraphics[width=0.48\textwidth]{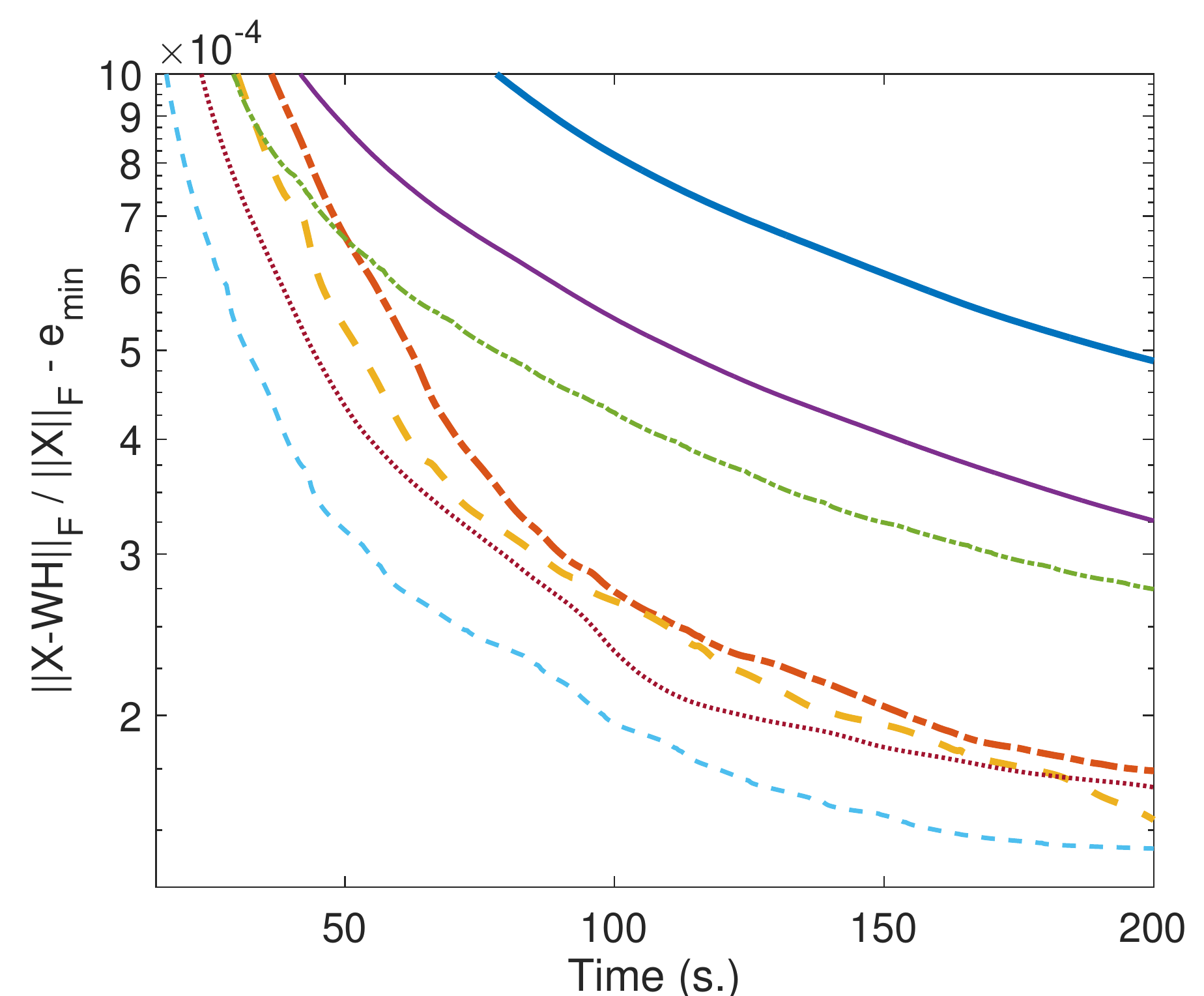}  & 
\includegraphics[width=0.48\textwidth]{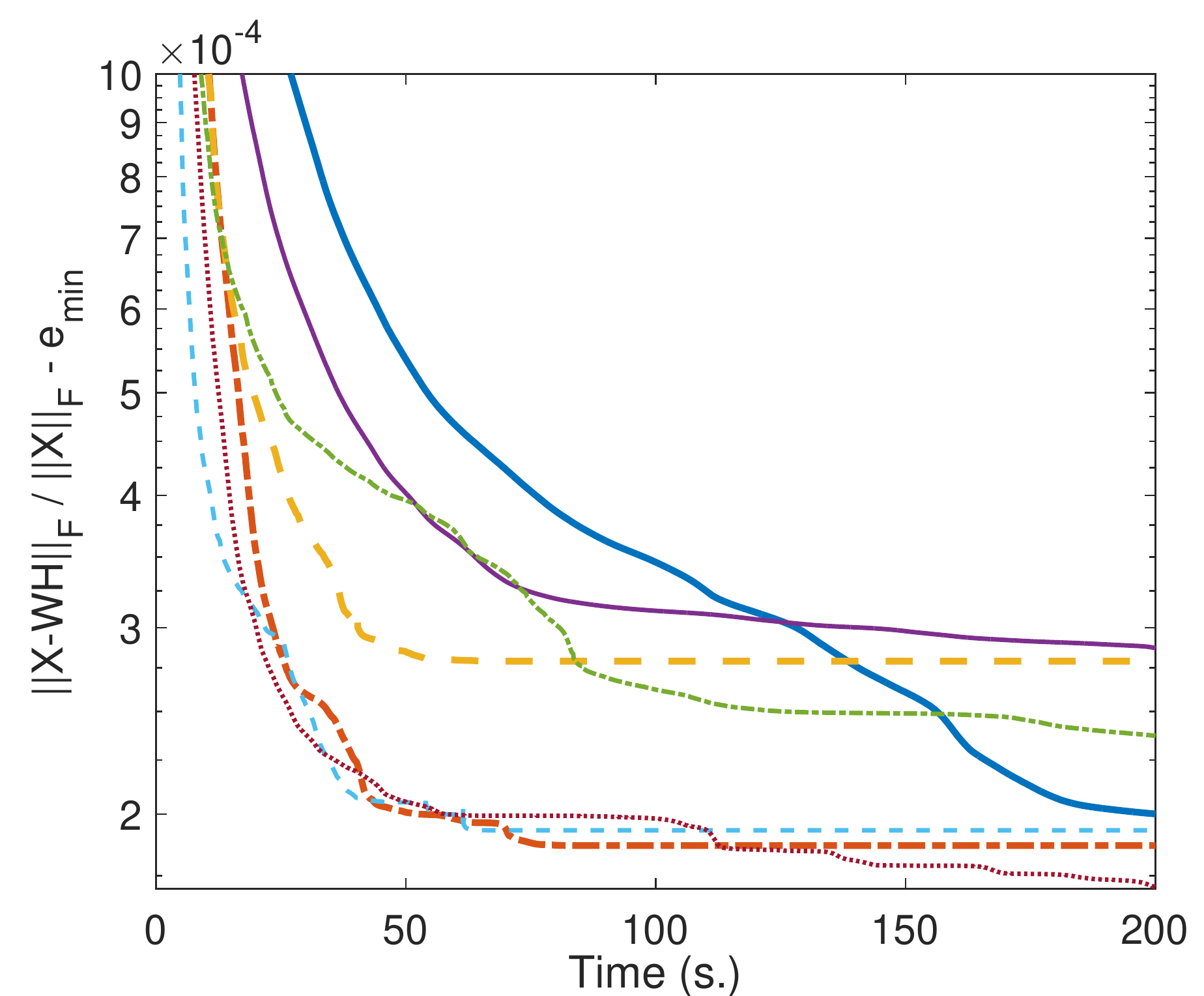} \\ 
\end{tabular}
\caption{
Average value of the error measure~\eqref{relerr} of ANLS, A-HALS and their extrapolated variants applied on  the 4 image data sets: 
CBCL (top left),
Umist (top right),  
ORL (bottom left), 
Frey (bottom right). 
   \label{fig:image}}
\end{center}
\end{figure}

We now run the algorithms for 200 seconds on the 4 image data sets; see 
Figure~\ref{fig:image} which displays the evolution of the average error measure~\eqref{relerr} for each data set, and Table~\ref{imdat} which compares the final errors obtained by the different algorithms.

\begin{center}  
 \begin{table}[h!] 
 \begin{center} 
\caption{Comparison of the final relative error obtained by the NMF algorithms on the image data sets: 
Average error, standard deviation and rankings among the 40 runs (4 data sets, 10 initializations each). 
The $i$th entry of the vector indicates the number of times the algorithm generated the $i$th best solution. \label{errim}} 
 \begin{tabular}{|c|c|c|c|c|c|c|c|} 
 \hline Algorithm &  mean $\pm$ std & ranking  \\ 
 \hline 
 ANLS &              $0.110703 \pm 2.964 \, 10^{-2}$ &  ( 3,  3,  3,  5,  5,  8, 13)   \\ 
E-ANLS ($hp=1$) &    
		     $\mathbf{0.110547 \pm 2.958\, 10^{-2}}$ &  ( 9, 12,  7,  5,  4,  2,  1)   \\ 
E-ANLS ($hp=3$) &    $0.110570 \pm 2.956\, 10^{-2}$ &  ( 9,  6,  5,  7,  2,  8,  3)   \\ 
 A-HALS &            $0.110690 \pm 2.956\, 10^{-2}$ &  ( 1,  4,  4,  2,  3, 13, 13)   \\ 
E-A-HALS ($hp=1$) &  $0.110634 \pm 2.958\, 10^{-2}$ &  ( 4,  2,  2,  4, 17,  7,  4)   \\ 
E-A-HALS ($hp=3$) &  $0.110552 \pm 2.956\, 10^{-2}$ &  ( 5, 10, 11,  8,  3,  0,  3)   \\ 
APG-MF &             $0.110559 \pm 2.956\, 10^{-2}$ &  ( 9,  3,  8,  9,  6,  2,  3)   \\ 
\hline \end{tabular} 
 \end{center} 
 \end{table} 
 \end{center} 
 
 We observe the following: 
 \begin{itemize}
 
 \item E-A-HALS ($hp=3$) has the fastest initial convergence speed, followed by E-ANLS variants and APG-MF. 
 As in the preliminary numerical experiments, E-A-HALS ($hp=1$) is able to accelerate A-HALS but not as much as E-A-HALS ($hp=3$). 

\item In terms of final error, there is no clear winner between the extrapolated variants (similarly as for the full-rank synthetic data sets), while ANLS that clearly performs the worst.  
 
 \end{itemize}

To conclude, we see that the extrapolation scheme is particularly beneficial to ANLS that is significantly accelerated and even able to outperform E-A-HALS in some cases (while A-HALS performs in general much better than ANLS, as already pointed out in~\cite{gillis2012accelerated}). 
Athough APG-MF outperforms ANLS (as already observed in~\cite{xu2013block}) and A-HALS, 
it is in general outperformed by the other extrapolated variants.

\subsubsection{ Sparse document data sets } \label{sec:text}

\begin{figure}[H]
\begin{center}
\begin{tabular}{cc}
\includegraphics[width=0.48\textwidth]{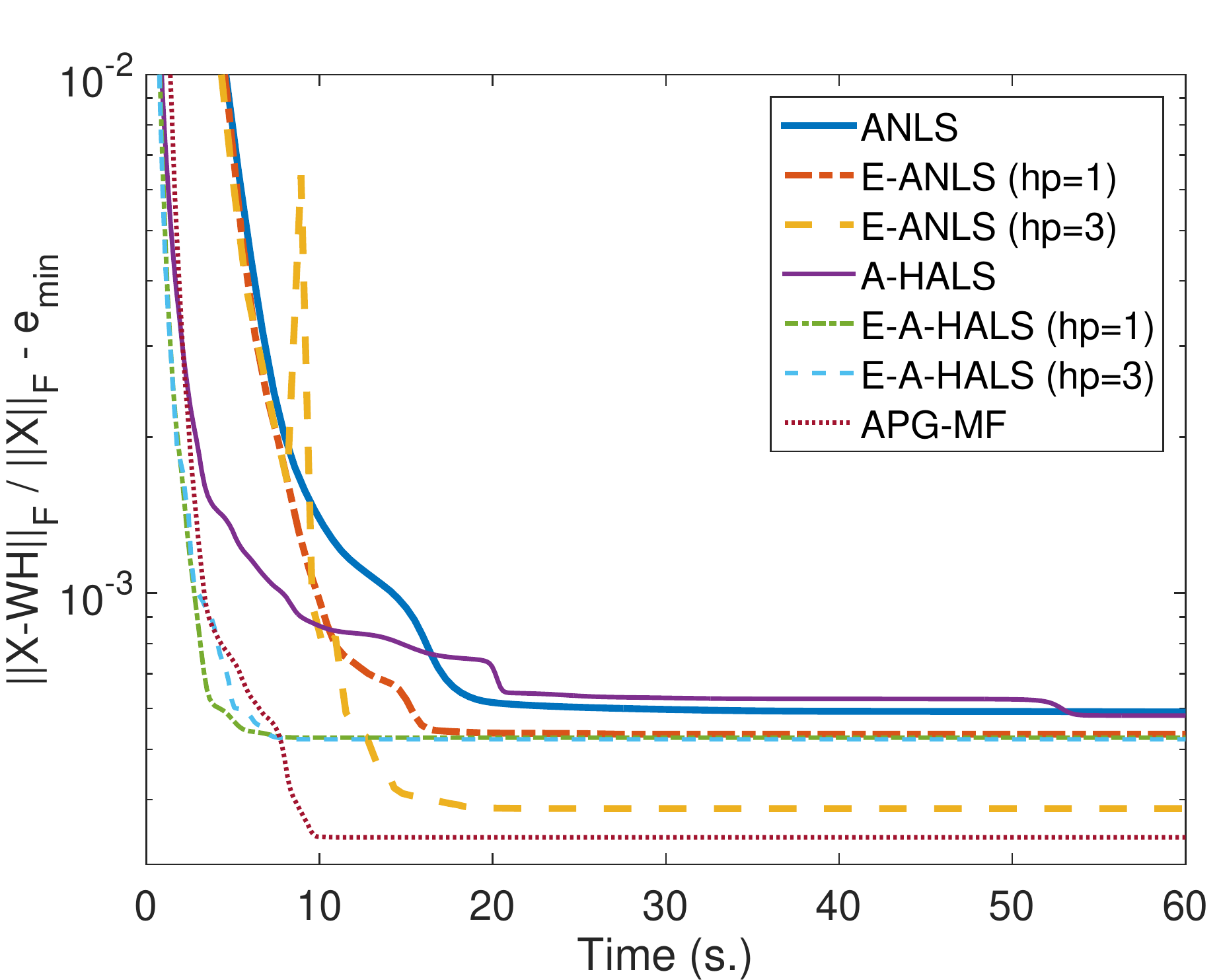}  & 
\includegraphics[width=0.48\textwidth]{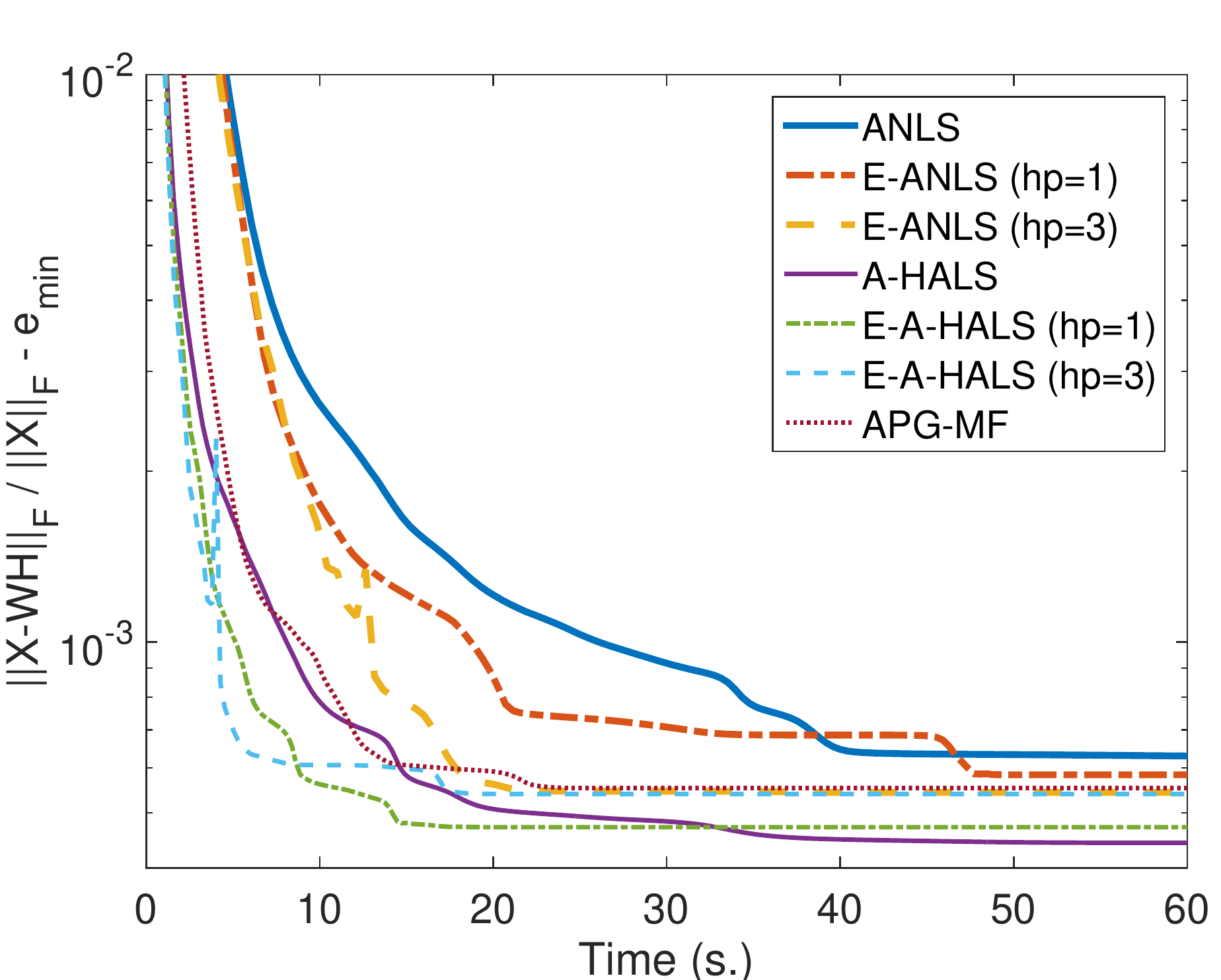} \\ 
\includegraphics[width=0.48\textwidth]{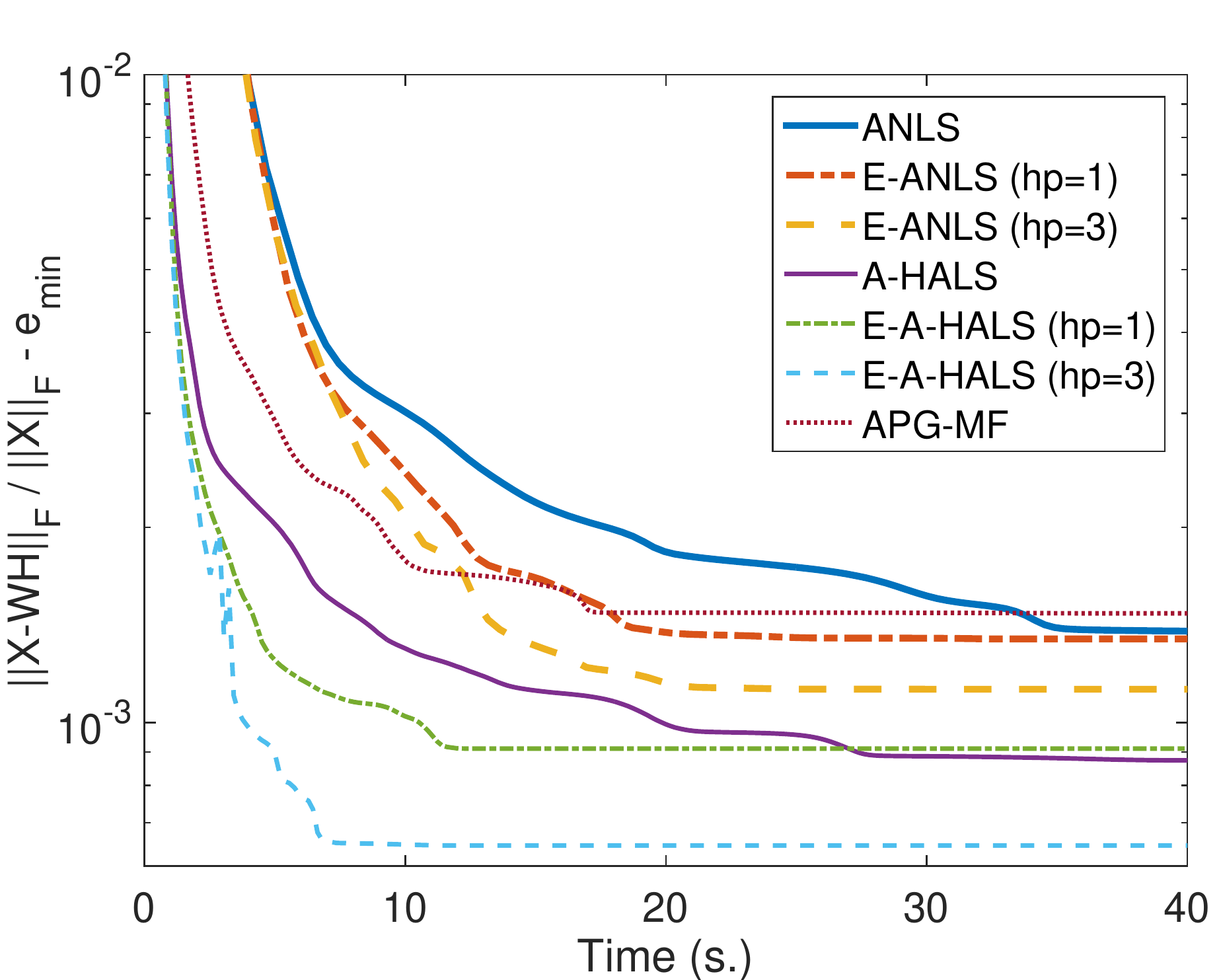}  & 
\includegraphics[width=0.48\textwidth]{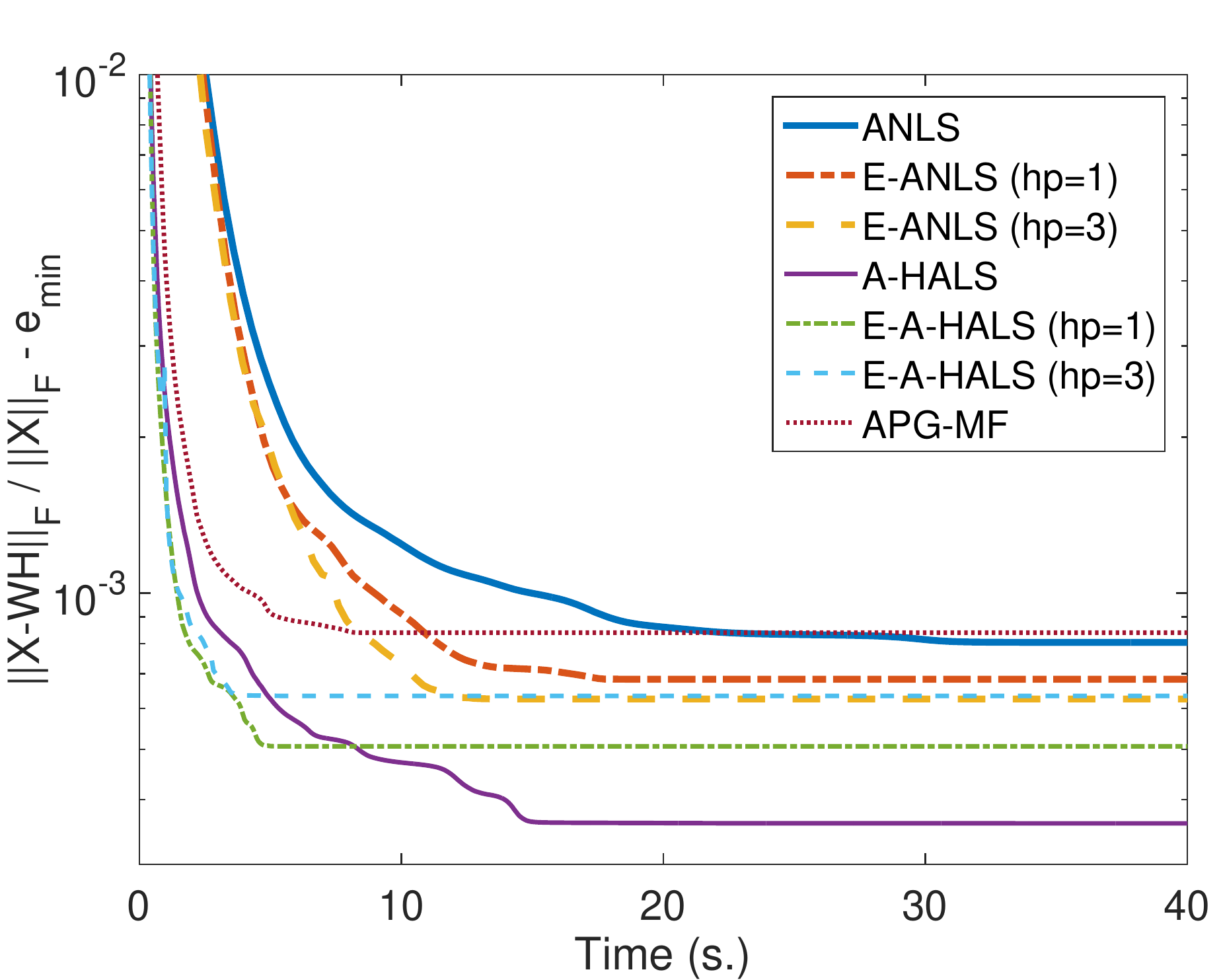} \\ 
\includegraphics[width=0.48\textwidth]{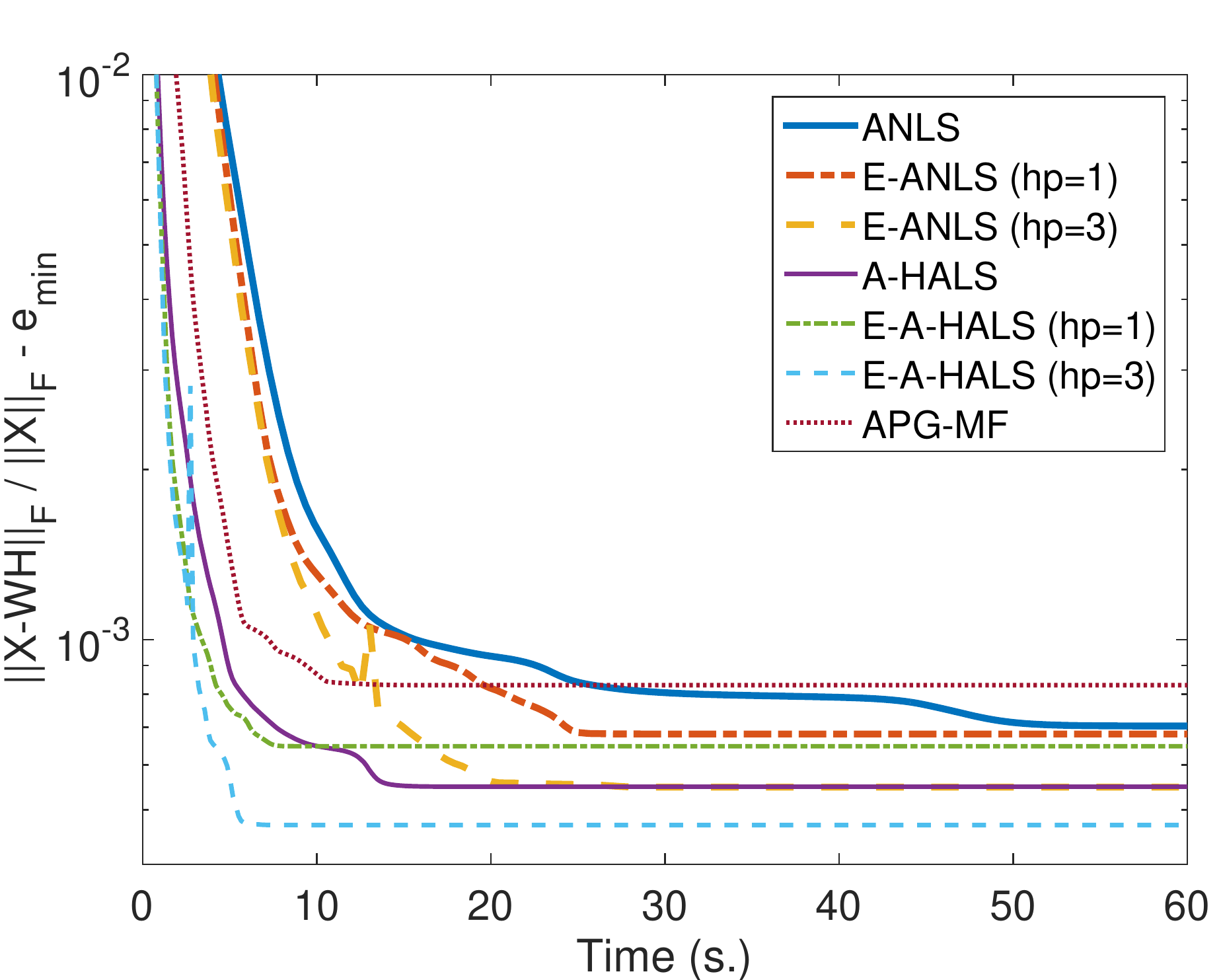}  & 
\includegraphics[width=0.48\textwidth]{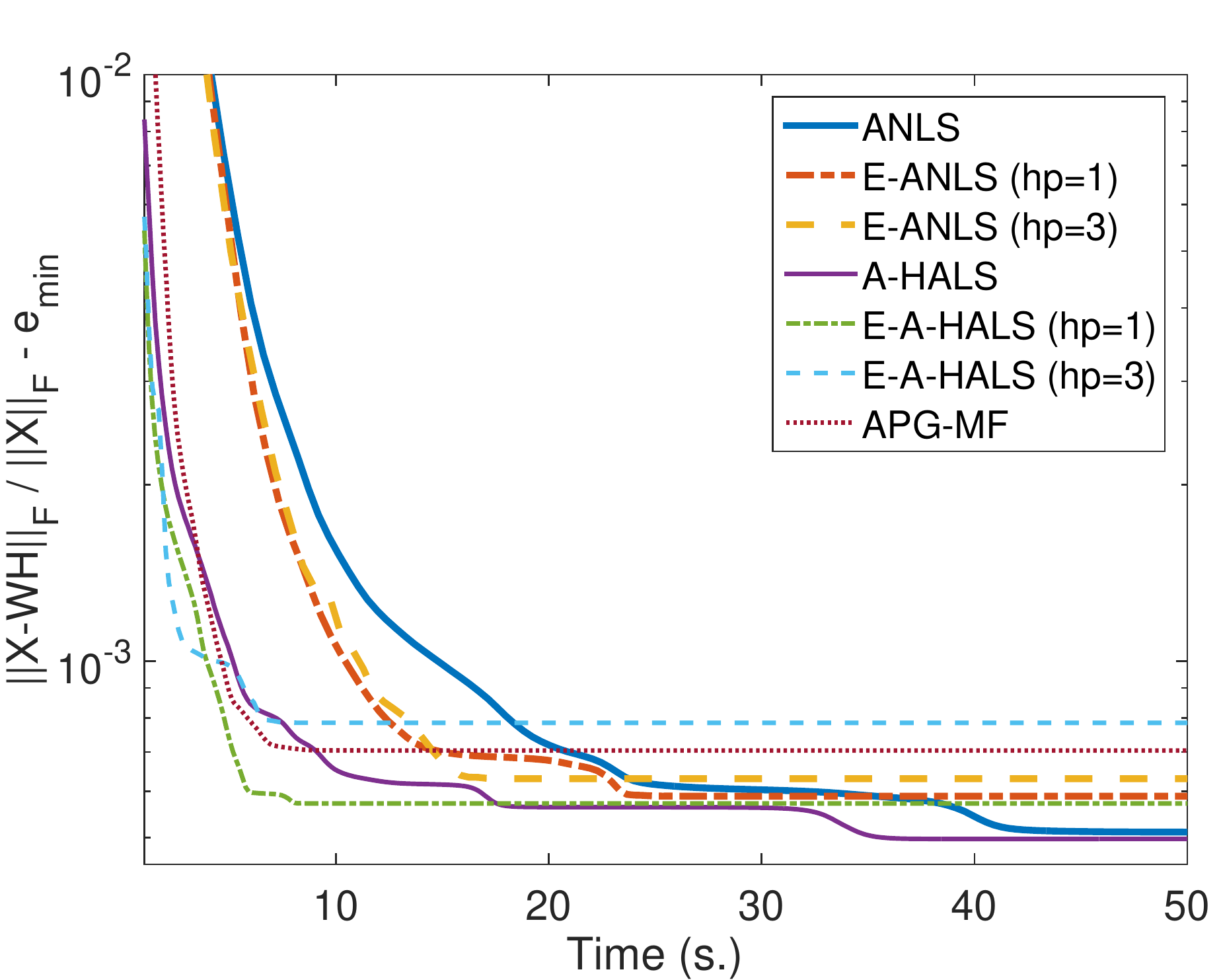} \\ 
\end{tabular}
\caption{
Average value of the error measure~\eqref{relerr} of ANLS, A-HALS and their extrapolated variants applied on  the 6 documents data sets: 
classic (top left), 
sports (top right), 
reviews (middle left), 
hitech (middle right), 
ohscal (bottom left), 
la1 (bottom right). 
  \label{fig:text}}
\end{center}
\end{figure}

We now run for 200 seconds ANLS, A-HALS and their extrapolated variants on the 6 document data sets; see 
Figure~\ref{fig:text} which displays the evolution of the average error measure~\eqref{relerr} for each data set, and Table~\ref{tab:text} which compares the final errors obtained by the different algorithms.

\begin{center}  
 \begin{table}[h!] 
 \begin{center} 
\caption{Comparison of the final relative error obtained by the NMF algorithms on the document data sets: 
Average error, standard deviation and rankings among the 60 runs (6 data sets, 10 initializations each). 
The $i$th entry of the vector indicates the number of times the algorithm generated the $i$th best solution. \label{tab:text} } 
 \begin{tabular}{|c|c|c|c|c|c|c|c|} 
 \hline Algorithm &  mean $\pm$ std & ranking  \\ 
 \hline 
 ANLS &              $0.850433 \pm 3.186 \, 10^{-2}$ &  ( 5,  3, 12,  6,  9, 11, 14)   \\ 
E-ANLS ($hp=1$) &    $0.850417 \pm 3.187\, 10^{-2}$ &  ( 7,  8,  6, 12, 13, 12,  2)   \\ 
E-ANLS ($hp=3$) &    $0.850324 \pm 3.189\, 10^{-2}$ &  ( 9, 11,  6,  9, 15,  6,  4)   \\ 
 A-HALS &            $\mathbf{0.850232 \pm 3.198\, 10^{-2}}$ &  (15, 11,  9,  8,  7,  7,  3)   \\ 
E-A-HALS ($hp=1$) &  $0.850287 \pm 3.198\, 10^{-2}$ &  (13, 13, 12,  6,  7,  6,  3)   \\ 
E-A-HALS ($hp=3$) &  $0.850281 \pm 3.204\, 10^{-2}$ &  (15, 11, 11,  4,  5,  9,  5)   \\ 
APG-MF &             $0.850471 \pm 3.183\, 10^{-2}$ &  ( 5,  5,  9, 10,  5,  9, 17)   \\ 
\hline \end{tabular} 
 \end{center} 
 \end{table} 
 \end{center}

 We observe the following: 
 \begin{itemize}
 
 \item E-A-HALS variants have the fastest initial convergence speed converging in average in about 10 seconds, followed by A-HALS which sometimes takes much more time to stabilize (e.g., more than 50 seconds for the classic data set).  
 APG-MF does not converge as fast as E-A-HALS variants. 
 E-ANLS variants converge much faster than ANLS but sometimes take more than 30 seconds to stabilize.

\item In terms of final error,  there is no clear winner although A-HALS and E-A-HALS ($hp=3$) gives most of the time the best solution (15 our of 60 cases). APG-MF tends to generate the worst solutions (17 out of the 60 cases) and performs similarly as ANLS in this respect.  
 
 \end{itemize}

For sparse data sets, E-A-HALS is the best option for which both variants ($hp=1,3$) perform similarly. APG-MF and ANLS and its extrapolated variants are less effective in this case.  

%To conclude, 
%we see that the extrapolation scheme is particularly beneficial to ANLS that is significantly accelerated and even %able to outperform E-A-HALS in some cases (while A-HALS performs in general much better than ANLS, as already pointed out in~\cite{gillis2012accelerated}). 

\section{Conclusion}

In this paper, we have proposed an extrapolation scheme for NMF algorithms to significantly accelerate their convergence.  We have focused on two state-of-the-art NMF algorithms, namely, ANLS~\cite{kim2011fast} and A-HALS~\cite{gillis2012accelerated}. The main conclusions are the following: 
\begin{itemize}

\item In all cases, the extrapolated variants significantly outperform the original algorithms. 

\item For randomly generated low-rank matrices, E-ANLS, the extrapolated variant of ANLS, allows a significant acceleration, being able to compute solutions with very small relative errors ($\approx 10^{-8}$) in all cases while the other approaches fail to do so. 

\item For dense data sets, E-ANLS and E-A-HALS perform similarly although A-HALS performs much better than ANLS. This is interesting: the extrapolated variants allowed ANLS to get back on A-HALS. 

\item For sparse data sets, E-A-HALS performs the best and should be preferred to the other variants. 

\item The extrapolated projected gradient method proposed in~\cite{xu2013block} and referred to as APG-MF performs well but does not perform as well as the extrapolated variants proposed in this paper.

\end{itemize}

This work was mostly experimental. 
It would be crucial to understand the extrapolation scheme better from a theoretical point of view. 
In particular, can we prove convergence to a stationary point as done in~\cite{xu2013block}? and can we quantify precisely the acceleration like it has been done in the convex case? 
Further work also include the use of extrapolation in other settings, e.g., 
NMF with other objective functions such as the $\beta$ divergence~\cite{fevotte2011algorithms}, 
nonnegative tensor factorization (NTF)~\cite{cichocki2009nonnegative}, and symmetric NMF~\cite{vandaele2016efficient}.

\bibliographystyle{spmpsci}  
\bibliography{FasterNMF}

\begin{thebibliography}{10}
\providecommand{\url}[1]{{#1}}
\providecommand{\urlprefix}{URL }
\expandafter\ifx\csname urlstyle\endcsname\relax
  \providecommand{\doi}[1]{DOI~\discretionary{}{}{}#1}\else
  \providecommand{\doi}{DOI~\discretionary{}{}{}\begingroup
  \urlstyle{rm}\Url}\fi

\bibitem{beck2013convergence}
Beck, A., Tetruashvili, L.: On the convergence of block coordinate descent type
  methods.
\newblock SIAM Journal on Optimization \textbf{23}(4), 2037--2060 (2013)

\bibitem{chambolle2015remark}
Chambolle, A., Pock, T.: A remark on accelerated block coordinate descent for
  computing the proximity operators of a sum of convex functions.
\newblock SMAI Journal of Computational Mathematics \textbf{1}, 29--54 (2015)

\bibitem{chow2017cyclic}
Chow, Y.T., Wu, T., Yin, W.: Cyclic coordinate-update algorithms for
  fixed-point problems: Analysis and applications.
\newblock SIAM Journal on Scientific Computing \textbf{39}(4), A1280--A1300
  (2017)

\bibitem{Cic4}
Cichocki, A., Phan, A.H.: {Fast local algorithms for large scale Nonnegative
  Matrix and Tensor Factorizations}.
\newblock IEICE Transactions on Fundamentals of Electronics \textbf{Vol. E92-A
  No.3}, 708--721 (2009)

\bibitem{Cic}
Cichocki, A., Zdunek, R., Amari, S.: {Hierarchical ALS Algorithms for
  Nonnegative Matrix and 3D Tensor Factorization}.
\newblock Lecture Notes in Computer Science, Springer \textbf{4666}, 169--176
  (2007)

\bibitem{cichocki2009nonnegative}
Cichocki, A., Zdunek, R., Phan, A.H., Amari, S.i.: Nonnegative matrix and
  tensor factorizations: applications to exploratory multi-way data analysis
  and blind source separation.
\newblock John Wiley \& Sons (2009)

\bibitem{erichson2018randomized}
Erichson, N.B., Mendible, A., Wihlborn, S., Kutz, J.N.: Randomized nonnegative
  matrix factorization.
\newblock Pattern Recognition Letters  (2018).
\newblock Doi:10.1016/j.patrec.2018.01.007

\bibitem{fercoq2015accelerated}
Fercoq, O., Richt{\'a}rik, P.: Accelerated, parallel, and proximal coordinate
  descent.
\newblock SIAM Journal on Optimization \textbf{25}(4), 1997--2023 (2015)

\bibitem{fevotte2011algorithms}
F{\'e}votte, C., Idier, J.: Algorithms for nonnegative matrix factorization
  with the $\beta$-divergence.
\newblock Neural computation \textbf{23}(9), 2421--2456 (2011)

\bibitem{fu2018nonnegative}
Fu, X., Huang, K., Sidiropoulos, N.D., Ma, W.K.: Nonnegative matrix
  factorization for signal and data analytics: Identifiability, algorithms, and
  applications.
\newblock arXiv preprint arXiv:1803.01257  (2018)

\bibitem{gillis2014}
Gillis, N.: {The Why and How of Nonnegative Matrix Factorization}.
\newblock In: J.~Suykens, M.~Signoretto, A.~Argyriou (eds.) Regularization,
  Optimization, Kernels, and Support Vector Machines, pp. 257--291. Chapman \&
  Hall/CRC, Machine Learning and Pattern Recognition Series (2014)

\bibitem{gillis2017}
Gillis, N.: Introduction to nonnegative matrix factorization.
\newblock SIAG/OPT Views and News \textbf{25}(1), 7--16 (2017)

\bibitem{gillis2012accelerated}
Gillis, N., Glineur, F.: Accelerated multiplicative updates and hierarchical
  {ALS} algorithms for nonnegative matrix factorization.
\newblock Neural computation \textbf{24}(4), 1085--1105 (2012)

\bibitem{nenmf}
Guan, N., Tao, D., Luo, Z., Yuan, B.: Nenmf: An optimal gradient method for
  nonnegative matrix factorization.
\newblock IEEE Transactions on Signal Processing \textbf{60}(6), 2882--2898
  (2012)

\bibitem{hong2017iteration}
Hong, M., Wang, X., Razaviyayn, M., Luo, Z.Q.: Iteration complexity analysis of
  block coordinate descent methods.
\newblock Mathematical Programming \textbf{163}(1-2), 85--114 (2017)

\bibitem{hsieh2011fast}
Hsieh, C.J., Dhillon, I.S.: Fast coordinate descent methods with variable
  selection for non-negative matrix factorization.
\newblock In: Proceedings of the 17th ACM SIGKDD international conference on
  Knowledge discovery and data mining, pp. 1064--1072. ACM (2011)

\bibitem{kim2008nonnegative}
Kim, H., Park, H.: Nonnegative matrix factorization based on alternating
  nonnegativity constrained least squares and active set method.
\newblock SIAM Journal on Matrix Analysis and Applications \textbf{30}(2),
  713--730 (2008)

\bibitem{kim2014algorithms}
Kim, J., He, Y., Park, H.: Algorithms for nonnegative matrix and tensor
  factorizations: A unified view based on block coordinate descent framework.
\newblock Journal of Global Optimization \textbf{58}(2), 285--319 (2014)

\bibitem{kim2011fast}
Kim, J., Park, H.: Fast nonnegative matrix factorization: An active-set-like
  method and comparisons.
\newblock SIAM Journal on Scientific Computing \textbf{33}(6), 3261--3281
  (2011)

\bibitem{lee1999learning}
Lee, D.D., Seung, H.S.: Learning the parts of objects by non-negative matrix
  factorization.
\newblock Nature \textbf{401}(6755), 788 (1999)

\bibitem{lin2007projected}
Lin, C.J.: Projected gradient methods for nonnegative matrix factorization.
\newblock Neural computation \textbf{19}(10), 2756--2779 (2007)

\bibitem{luenberger1984linear}
Luenberger, D.G., Ye, Y.: Linear and nonlinear programming, Fourth edition.
\newblock Springer (2015).
\newblock Available
  from~\url{https://web.stanford.edu/class/msande310/310trialtext.pdf}

\bibitem{nesterov2013introductory}
Nesterov, Y.: Introductory lectures on convex optimization: A basic course,
  vol.~87.
\newblock Springer Science \& Business Media (2013)

\bibitem{o2015adaptive}
O’donoghue, B., Cand\`es, E.: Adaptive restart for accelerated gradient
  schemes.
\newblock Foundations of computational mathematics \textbf{15}(3), 715--732
  (2015)

\bibitem{o2017behavior}
O'Neill, M., Wright, S.J.: Behavior of accelerated gradient methods near
  critical points of nonconvex problems.
\newblock arXiv preprint arXiv:1706.07993  (2017)

\bibitem{paquette18a}
Paquette, C., Lin, H., Drusvyatskiy, D., Mairal, J., Harchaoui, Z.: Catalyst
  for gradient-based nonconvex optimization.
\newblock In: A.~Storkey, F.~Perez-Cruz (eds.) Proceedings of the Twenty-First
  International Conference on Artificial Intelligence and Statistics,
  \emph{Proceedings of Machine Learning Research}, vol.~84, pp. 613--622. PMLR,
  Playa Blanca, Lanzarote, Canary Islands (2018)

\bibitem{vandaele2016efficient}
Vandaele, A., Gillis, N., Lei, Q., Zhong, K., Dhillon, I.: Efficient and
  non-convex coordinate descent for symmetric nonnegative matrix factorization.
\newblock IEEE Transactions on Signal Processing \textbf{64}(21), 5571--5584
  (2016)

\bibitem{vavasis2009complexity}
Vavasis, S.A.: On the complexity of nonnegative matrix factorization.
\newblock SIAM Journal on Optimization \textbf{20}(3), 1364--1377 (2010)

\bibitem{xu2013block}
Xu, Y., Yin, W.: A block coordinate descent method for regularized multiconvex
  optimization with applications to nonnegative tensor factorization and
  completion.
\newblock SIAM Journal on Imaging Sciences \textbf{6}(3), 1758--1789 (2013)

\bibitem{ZG05}
Zhong, S., Ghosh, J.: Generative model-based document clustering: a comparative
  study.
\newblock Knowledge and Information Systems \textbf{8 (3)}, 374--384 (2005)

\end{thebibliography}

\end{document}